\documentclass[11pt]{article}

\usepackage[margin=1in]{geometry} 
\usepackage{amsmath,amssymb,amsfonts,amsthm}
\usepackage{graphicx}
\usepackage{float}
\usepackage{xcolor}
\usepackage[authoryear, round]{natbib}
\usepackage{hyperref}
\hypersetup{colorlinks,citecolor=blue,linkcolor=blue,breaklinks=true}
\usepackage{epstopdf}
\usepackage{yhmath} 
\usepackage{makecell}
\usepackage{enumerate}
\usepackage{enumitem}
\usepackage{natbib}

\allowdisplaybreaks

\newtheorem{theorem}{Theorem}[section]
\newtheorem{lemma}{Lemma}[section]

\newtheorem{proposition}[theorem]{Proposition}
\newtheorem{definition}{Definition}[section]

\newtheorem{assumption}{\textbf{Assumption}}[section]
\newcommand{\R}{\mathbb{R}}

\newcommand{\bbm}{\begin{bmatrix}}
	\newcommand{\ebm}{\end{bmatrix}}

\numberwithin{equation}{section}

\begin{document}
	
	\title{Mean field game problem for the optimal control of neuronal spiking activity}
 
	\author{
		Lijun Bo\thanks{Email: lijunbo@xidian.edu.cn,  School of Mathematics and Statistics, Xidian University, Xi'an, 710126, China.}
		\and
		Dongfang Yang\thanks{Email: yangdf@stu.xidian.edu.cn,  School of Mathematics and Statistics, Xidian University, Xi'an, 710126, China.}
  \and
		Shihua Wang\thanks{Corresponding author. Email: wangshihua@xidian.edu.cn,  School of Mathematics and Statistics, Xidian University, Xi'an, 710126, China.}
	}
	
	\date{}
	
	\maketitle

    \begin{abstract}
    	We study the mean field game problem for a nervous system consisting of a large number of neurons with mean-field interaction. In this system, each neuron can modulate its spiking activity by controlling its membrane potential to synchronize with others, thereby giving rise to a finite-player game problem. To address this, we first examine the corresponding mean field game problem and characterize the mean field equilibrium by solving a fixed point problem. Subsequently, leveraging the obtained mean field equilibrium, we construct an approximate Nash equilibrium for the finite-player game as the number of neurons is large.
     
\vspace{0.1in}

\vspace{0.02in}
    \end{abstract}

\section{Introduction}

Within biological nervous systems, a large number of neurons exhibit intricate behaviors and interactions (c.f. \cite{G02}, \cite{H09}, \cite{delarue2015particle}, and \cite{grazieschi2019}). In recent decades, extensive research has been dedicated to unraveling the structures and mechanisms underlying these complex systems. Various models have been designed to emulate the characteristic behaviors of neurons.  In the context of individual neurons, widely employed models encompass the Hodgkin-Huxley model in \cite{HH52} and the Leaky Integrate-and-Fire (LIF) model in \cite{L07}, along with their various adaptations. Considering the discharge characteristics, there is a long-standing tradition of utilizing finite systems of point processes in both discrete and continuous time to model neuronal dynamics. Building upon this foundation, researchers have developed more generalized processes for modeling purposes. \cite{C11} examine finite-state stochastic chains within the context of the leaky integrate-and-fire model. \cite{GL13} expand on these interacting chains by incorporating memory of variable length. Furthermore, \cite{delarue2015global} investigate the well-posedness of a networked integrate-and-fire model describing an infinite population of neurons which interact with each other through their common statistical distribution.
The relevant literature also includes \cite{DG15}, \cite{GL16}, and \cite{FL16}, among others.
    
The regulation of neuronal activity constitutes another pivotal area within the realms of neuroscience research and medical intervention.  In the work outlined by \cite{FG07}, addressing Parkinson’s disease entails modulating the firing rates of neurons within targeted brain areas. 
As technology progresses, numerous techniques for modulating neuronal activity have been developed. These include direct electrical stimulation (e.g., \cite{SB90}, \cite{DM03}), optical methods (e.g., \cite{CY02}, \cite{BZ05}), and micro-injection techniques to inhibit neuronal activity (e.g., \cite{DZ01}). The field has already extensively explored the resulting optimal control problem. For example, \cite{FT03} investigate an optimal control issue within a LIF neuron model, with the aim of minimizing the variance of the membrane potential at spiking times, all while ensuring the mean remains constant. Similarly, \cite{AP11} seek to determine the optimal external control of the membrane potential of an LIF neuron to achieve precise spike timings. The optimal control problems for single neurons or coupled neurons have also been thoroughly examined in \cite{LD13}, \cite{LW19}, \cite{MG22}, and the references therein.
    
 Since the seminal work of \cite{HM06} and \cite{LL07}, which was conducted independently, the mean field game (MFG) theory has primarily focused on examining games within the context of large populations. For more comprehensive insights into MFG theory, the reader may refer to \cite{BS16}, \cite{CD18}, and the additional references cited within those works. On the other hand, in recent years, there has been a growing interest in MFGs with the underlying discontinuous state processes.
 \cite{HA14} explore stochastic control problems involving nonlinear controlled jump diffusion processes of mean-field type, also known as McKean-Vlasov equations. \cite{BC19, BC20} consider MFGs with controlled jump diffusion processes and provide an approximate Nash equilibrium for the $n$-player game. Additionally, \cite{CF20} employs a probabilistic representation approach to study MFGs and the corresponding $n$-player game with finite states, where the system dynamics in terms of stochastic differential equations are driven by Poisson random measures. \cite{BL22} study a class of dynamic games with finite agents under a stochastic growth model, in which each agent's capital stock dynamics are characterized by a jump process. \cite{BW24} handle MFG problems with controlled jumps and common noise in the context of financial markets, where the jump risk is modeled by a self-exciting Hawkes process. 

In this paper, we study the optimal external control problem in a large population nervous system, in which each neuron interacts with others through both chemical and electrical synapses in a mean field fashion. As in \cite{GL13, GL16}, neurons spike at random times, causing a sharp decrease in their membrane potential. Concurrently, chemical synapses increase the potential of neighboring neurons, while electrical synapses adjust the membrane potential towards the system’s average value via gap junctions at a constant speed. From a mathematical perspective, the state of the system, representing the membrane potentials of the neurons, is characterized by controlled jump-diffusion processes that are interconnected through their drift terms and the deriving Poisson random measures. Since normal synchronization within neural populations is crucial for effective information transmission (c.f. \cite{A06}, \cite{WF07} and \cite{AP94}), our goal for each neuron is to fine-tune its external control in order to align its membrane potential as closely as possible with the system's average value. This gives rise to a finite-player stochastic game, which we tackle using the MFG approach. The solution procedure for the MFG is bifurcated into two distinct stages: initially, we resolve a conventional stochastic control problem as the number of neurons tends to infinity; subsequently, we ascertain the existence and uniqueness of a fixed point that pertains to the consistency condition. Having delineated the mean field equilibrium, we proceed to devise an approximate Nash equilibrium tailored for the finite-player game scenario. The proposed model in the paper, grounded in the essential characteristics of neurons, strikes an optimal balance between authenticity and problem tractability. By incorporating connection strength function into the model, we introduce a layer of complexity that more accurately reflects the intricate dynamics of neuronal interactions. However, this added complexity renders the problem more challenging to solve, particularly evident in establishing consistency conditions for mean field equilibrium. In this study, we deal with these consistency conditions, formulating them as a fixed point system. In particular, in our linear connection strength function case as in \cite{cormier2020long}, this complex system simplifies to a single one. Our approach ingeniously marries the contraction mapping principle with the implicit function theorem, thereby rigorously proving the existence and uniqueness of fixed points in the $C^1$ sense. A key advantage of our methodology lies in the explicit formulation of the best response solution, enabling us to efficiently determine the mean field equilibrium. This paves the way for constructing an approximate Nash equilibrium, significantly advancing our understanding and manipulation of neural network behaviors in complex systems.

This paper is structured as follows: In Section \ref{sec:model}, we formulate the controlled jump-diffusion dynamics of the membrane potential for interacting neurons and formulates the finite $n$-neuron game. In Section \ref{sec:MFE}, we establish the mean field equilibrium (MFE) in analytical form through the study of a fixed point problem by using the MFG approach. In Section~\ref{sec:approximateNE}, we construct and verify an approximate Nash equilibrium in terms of the obtaining MFE. The proof of the auxiliary result is reported in Appendix~\ref{sec:appendix}.

\section{The Model}\label{sec:model}
    
In this section, we describe the dynamics of membrane potential for interacting neurons as a controlled jump-diffusion system and formulate the corresponding finite $n$-player game.  

Let $(\Omega, \mathcal{F}, \mathbb{F}, \mathbb{P})$ be a complete filtered probability space with the filtration $\mathbb{F}=(\mathcal{F}_t)_{t\in [0,T]}$ satisfying the usual conditions, which supports independent Poisson random measures $(N^i(dt,dz))_{i=1}^n$ with a common intensity measure $\nu(dz)dt$ on $[0,T]\times{\cal B}([0,1])$. The jump measure $\nu(dz)$ satisfies $\int_0^1 z\nu(dz)<\infty$. Denote by $U_t^i$ the membrane potential of neuron $i$ at time $t\in[0,T]$ with $T>0$ being the finite terminal horizon. The dynamics of the membrane potential $\boldsymbol{U}=(U_t^1,\ldots, U_t^n)_{t\in[0, T]}$ is described as the following controlled jump-diffusion system, for $t\in[0,T]$,
    \begin{align}\label{ndynamics}
    U_t^i&=u_i-a_i\int_0^t\left(U_s^i-\frac{1}{n}\sum_{j=1}^nU_s^j\right)ds+\int_0^tc_i\theta_s^ids-\int_0^t\int_0^1zU_{s-}^iN^i(ds,dz)\nonumber\\&\quad+\frac{1}{n}\sum_{j\neq i}\int_0^t\int_0^1 \phi(zU_{s-}^j)N^j(ds,dz),
    \end{align}
where $U_0^i=u_i\in \mathbb{R}$ is the initial membrane potential value of neuron $i$, $a_i>0$ is the parameter involved with the gap junction, which brings a drift of the membrane potential of neuron $i$ to the average value, and similar to \cite{cormier2020long}, $\phi:\mathbb{R}\to \mathbb{R}$ is the connection strength between two interacting neurons, which is given by the linear mapping $\phi(x)=kx+\ell$ with parameters $k,\ell>0$.
    
The last two terms in SDE~\eqref{ndynamics} reflect the effect of spike activities in the system. At its spiking time, the membrane potential of the spiking neuron $i$ jumps. Simultaneously, it receives an additional amount of potential from other neurons interacting with it by chemical synapses. This amount depends on the connection strength function. We omit $n$ in the superscript for simplicity, which usually represents the total number of neurons in the system. Note that the membrane potential of the $i$-th neuron is affected by the others through the average value $\frac{1}{n}\sum_{j=1}^n U_s^j$ and the additional income. Here, the $\mathbb{F}$-adapted process $\boldsymbol{\theta}=(\theta_t^1,\ldots,\theta_t^n)_{t\in[0,T]}\in\mathcal{A}^n$ represents the external control to regulate the membrane potential of neurons, while $c_i$ is the scale factor of the control $\theta^i$ for neuron $i$, which determines the extent to which the implementation of external stimuli results in change. In particular, when the scale factor $c_i=0$ (i.e., without control), the model \eqref{ndynamics} is reduced to that in \cite{GL13} and \cite{GL16}. The set $\mathcal{A}$ is the admissible control space, which is the collection of the c\`{a}dl\`{a}g and $\mathbb{F}$-adapted real-valued process $\theta= (\theta_t)_{t\in[0,T]}$ satisfying $\|\theta\|_{2,T}:=\sqrt{\mathbb{E}[\int_0^T\vert\theta_t\vert^2dt]}<+\infty$.

Next, we introduce the objective functionals for $n$ neurons in the interacting system by extending the ones in the existing literature, which cover the minimization of the variance of the membrane potential as in \cite{FT03} and the control cost as in \cite{ID14}. For any $i=1,\ldots, n$, we aim to minimize the cost functional for neuron $i$ in the system as follows: 
\begin{align}\label{objective}
    \inf_{\theta^i\in \mathcal{A}}J_i(\boldsymbol{\theta}):=\inf_{\theta^i\in \mathcal{A}}\mathbb{E}\left[\int_0^Tf^i(\boldsymbol{U}_t, \theta^i_t)dt+g^i(\boldsymbol{U}_T)\right],
\end{align}
where the running cost function $f^i$ and the terminal cost function $g^i$ are respectively given by $f^i(x,\theta):=\theta^2+\rho\theta(x^i-\bar{x})+\beta(x^i-\bar{x})^2$ and $g^i(x):=\gamma(x^i-\bar{x})^2$ with $\bar{x}=\frac{1}{n}\sum_{j=1}^nx^j$. Here, $\beta,\gamma>0$ are parameters which balance the membrane potential of this neuron with the average value of all neurons in the system, $\rho>0$ is the parameter that weighs the contributions of the relative sizes of the membrane potential on the choice of the control $\theta$. Here, we assume $\rho^2<4\beta$ in order to guarantee the convexity of the running cost function.
	
Due to the presence of the term $\bar{U}_t:=\frac{1}{n}\sum_{j=1}^n U_t^j$ in \eqref{ndynamics} and \eqref{objective}, solving the $n$-player game and obtaining the exact Nash equilibrium in closed form is unfeasible. Instead, we study the approximate Nash equilibrium, which is defined as follows.
\begin{definition}[Approximate Nash Equilibrium]
Let $\boldsymbol{\theta}^{-i}:=(\theta^1,\ldots ,\theta^{i-1},\theta^{i+1},\ldots,\theta^n)$. An admissible strategy $\boldsymbol{\theta}^{*}=(\theta^{*,1},\ldots \theta^{*,n})\in \mathcal{A}^n$ is called an $\epsilon_n$-Nash equilibrium for the $n$-player game problem \eqref{ndynamics}-\eqref{objective} if, for any $(\theta^i)_{i=1}^n\in \mathcal{A}$, it holds that 
\begin{equation}
    J_i(\boldsymbol{\theta}^*)-\epsilon_n\leq \inf_{\theta^i\in \mathcal{A}}J_i\left(\theta^i, \boldsymbol{\theta}^{*,-i}\right), \quad \forall i=1, \ldots, n,
    \end{equation}
 where the error term $\epsilon_n$ satisfies  $\lim_{n\to\infty}\epsilon_n=0$.  
\end{definition}
    
In order to obtain the approximate Nash equilibrium, we collect the type vector of the model $p_i:=(u_i, a_i,c_i)$ for $i=1,\ldots,n$, and these type vectors induce an empirical measure given by $\mu_n:=\frac{1}{n}\sum_{i=1}^n\delta_{p_i}$ on $\mathcal{B}(\mathcal{O})$. Here, $\delta_x$ denotes the Dirac measure at the point $x$, and the set $\mathcal{O}$, called the type space, is a compact subset in $(0, +\infty)^3$. We then impose the following assumption:
\begin{assumption}\label{Ap}
There exists $\mu\in\mathcal{P}_2(\mathcal{O})$ such that the empirical measure $\mu_n$ weakly converges to $\mu$ as $n\to\infty$, that is, $\int_{\mathcal{O}}\varphi d\mu_n\to \int_{\mathcal{O}}\varphi d\mu $ as $n\to\infty$, for all $\varphi\in C_b(\mathcal{O})$. Here, $\mathcal{P}_2(\mathcal{O})$ denotes the set of probability measures on $\mathcal{O}$ with finite second-order moments.
    \end{assumption} 
We will study the MFG problem and characterize the mean-field equilibrium under Assumption \ref{Ap} in the forthcoming section. 
    
\section{Mean Field Equilibrium}\label{sec:MFE} 

In this section, we formulate the limiting mean field game problem and seek the mean field equilibrium as the number of neurons goes to infinity.
    
To formulate the MFG, let $N(dt,dz)$ be a Poisson random measure on $[0,T]\times[0,1]$ with intensity measure $\nu(dz)dt$ satisfying $\int_0^1z\nu(dz)<+\infty$, and $p=(u,a,c)$ be a deterministic sample from the probability distribution $\mu$ given in Assumption \ref{Ap}. We consider fixing $m_U^{(p)}=(m_U^{(p)}(t))_{t\in[0,T]}$ and $m_{\phi}^{(p)}=(m_{\phi}^{(p)}(t))_{t\in[0,T]}$ which are in $\mathcal{C}_T^1:=C^1([0,T];\R)$. Thus, for the given type vector $p=(u,a,c)\in\mathcal{O}$, we introduce the dynamics of the membrane potential of a representative neuron as follows:
\begin{equation}\label{dynamic}
	U_t^{(p)}=u-a\int_0^t\left(U_s^{(p)}-m_U^{(p)}(s)\right)ds+\int_0^tc\theta_sds-\int_0^t\int_0^1zU_{s-}^{(p)}N(ds, dz)+\int_0^tm_{\phi}^{(p)}(s)ds.
\end{equation}
The stochastic control problem of the representative neuron is given by
\begin{align}\label{MFC}
\inf_{\theta \in \mathcal{A}}&\bar{J}(\theta; m_U^{(p)}, m_{\phi}^{(p)})\\
&\quad:=\mathbb{E}\left[\int_0^T\left(\theta_s^2+\rho\theta_s\left(U_s^{(p)}-m_U^{(p)}(s)\right)+\beta\left|U_s^{(p)}-m_U^{(p)}(s)\right|^2\right)ds+\gamma\left|U_T^{(p)}-m_U^{(p)}(T)\right|^2\right].\nonumber
\end{align}

We shall establish the mean field equilibrium by solving the stochastic control problem \eqref{MFC} and a fixed point problem related to the so-called 
consistency condition. The definition of conditional mean field equilibrium is provided below:
\begin{definition}[Conditional Mean Field Equilibrium]\label{def:MFE}
For $p \in \mathcal{O}$, given $m_U^{(p)},m_{\phi}^{(p)}\in\mathcal{C}_T^1$, let $\theta^{*,(p)} \in \mathcal{A}$ be the best response to control problem \eqref{dynamic}-\eqref{MFC}. The control strategy $\theta^{*,(p)}\in\mathcal{A}$ is called a conditional mean field equilibrium if it is the best response to itself in the sense that $m_U^{(p)}(t)=\mathbb{E}[U_t^{*,(p)}]$ and $m_{\phi}^{(p)}(t)=\mathbb{E}[\int_0^1 \phi(z U_t^{*,(p)})\nu(dz)]$ for all $t\in [0,T]$. Here, $U^{*,(p)}=(U_t^{*,(p)})_{t\in[0,T]}$ is the membrane potential under the optimal control $\theta^{*,(p)}$. 
\end{definition}

\subsection{The best response solution}

Let us employ the dynamic programming principle to tackle the stochastic control problem \eqref{dynamic}-\eqref{MFC}. For all $(t,x)\in [0,T]\times \mathbb{R}$ and $p\in \mathcal{O}$, we define the value function as follows: 
\begin{align*}
V^{(p)}(t,x):=\inf_{\theta\in\mathcal{A}}\mathbb{E}_{t, x}\left[\int_t^T\left(\theta_s^2+\rho\theta_s(U_s^{(p)}-m_U^{(p)}(s))+\beta\left|U_s^{(p)}-m_U^{(p)}(s)\right|^2\right)ds+\gamma\left|U_t^{(p)}-m_U^{(p)}(T)\right|^2\right] 
\end{align*}
with the conditional expectation $\mathbb{E}_{t,x}[\cdot]:=\mathbb{E}[\cdot \vert U_t^{(p)}=x]$. 
By dynamic programming principle, the associated HJB equation satisfied by the value function $V^{(p)}(t,x)$ is given by
    \begin{align}\label{HJB}
		0&=
		V_t^{(p)}+\inf_{\theta\in\mathcal{A}}\Big\{V_x^{(p)}(-a(x-m_U^{(p)}(t))+ m_{\phi}^{(p)}(t)+c\theta)+\theta^2+\rho\theta(x-m_U^{(p)}(t))+\beta(x-m_U^{(p)}(t))^2\nonumber\\
		&\quad+\int_0^1\left(V^{(p)}(t,x(1-z))-V^{(p)}(t,x)\right)\nu(dz)\bigg\}
    \end{align}
with the terminal condition $V^{(p)}(T, x)=\gamma(x-m_U^{(p)}(T))^2$ for all $x\in \mathbb{R}$. The following lemma gives the optimal value function and feedback control in the closed form. 
\begin{lemma}\label{lem:bestresponse}
The value function $V^{(p)}(t,x)$ admits the closed-form given by
\begin{align}\label{eq:valuefcn}
V^{(p)}(t,x)=A_p(t)\left(x-m_U^{(p)}(t)\right)^2+B_p(t, m_U^{(p)}, m_{\phi}^{(p)})\left(x-m_U^{(p)}(t)\right)+C_p(t, m_U^{(p)}, m_{\phi}^{(p)}), 
\end{align}
and the optimal	(feedback) control function is given by
\begin{equation}\label{optimalcontrol}
\theta^{*,(p)}(t,x)=\left(-cA_p(t)-\frac{\rho}{2}\right)(x-m_U^{(p)}(t))-\frac{c}{2}B_p(t, m_U^{(p)}, m_{\phi}^{(p)}).
\end{equation}
Here, the coefficients are given by
\begin{equation}\label{solution}
\left\{
\begin{aligned}\displaystyle &A_p(t)=\frac{(\frac{\rho^2}{4}-\beta)(e^{(\delta^+-\delta^-)(T-t)}-1)-\gamma( \delta^+e^{(\delta^+-\delta^-)(T-t)})-\delta^-)}{\delta^-e^{(\delta^+-\delta^-)(T-t)}-\delta^+-c^2\gamma(e^{(\delta^+-\delta^-)(T-t)}-1)},\\
&B_p(t, m_U^{(p)}, m_{\phi}^{(p)})=2e^{-\int_t^T(a+c^2A_p(s)+\int_0^1z\nu(dz)+\frac{\rho c}{2})ds}\Bigg[\int_t^TA_p(s)e^{\int_s^T(a+c^2A_p(v)+\int_0^1z\nu(dz)+\frac{\rho c}{2})dv}\\
&\qquad\qquad\qquad\qquad \times\left(-(m_U^{(p)})'(s)+m_{\phi}^{(p)}(s)+m_U^{(p)}(s)\int_0^1(z^2-z)\nu(dz)\right)ds\Bigg],\\
&C_p(t, m_U^{(p)}, m_{\phi}^{(p)})=
\int_t^T\bigg\{ \int_0^1z^2\nu(dz)A_p(s)(m_U^{(p)})^2(s)+2A_p(s)m_{\phi}^{(p)}(s)m_U^{(p)}(s)\\
&\qquad\qquad\qquad\qquad  -\frac{1}{4}c^2B_p^2(s, m_U^{(p)}, m_{\phi}^{(p)}) -B_p(s, m_U^{(p)}, m_{\phi}^{(p)})(m_U^{(p)})'(s)\bigg\}ds,
\end{aligned}\right.
\end{equation}
where the constants 
            \begin{align}\label{para}
            \begin{cases}
                \displaystyle R:=\left(a+\frac{\rho c}{2}-\frac{1}{2}\int_0^1(z^2-2z)\nu(dz)\right)^2-c^2\left(\frac{1}{4}\rho^2-\beta\right)>0,\\[0.8em]
                \displaystyle \delta^{\pm}:=-\left(a+\frac{\rho c}{2}-\frac{1}{2}\int_0^1(z^2-2z)\nu(dz)\right)\pm\sqrt{R}.
            \end{cases}
	    \end{align}
    \end{lemma}
    
\begin{proof}
Firstly, assume that the HJB equation \eqref{HJB} has a classical solution which is strictly convex  w.r.t. $x$. Then, it follows from the first-order condition that, the optimal (feedback) control is expressed as: 
    \begin{equation}\label{FOC}
		\theta^{*,(p)}(t,x)=-\frac{c}{2}V_x^{(p)}(t,x)-\frac{\rho}{2}(x-m_U^{(p)}(t)).
    \end{equation}
By substituting \eqref{FOC} into Eq.~\eqref{HJB}, it holds that
    \begin{align}\label{pde}
		0&=
		V_t^{(p)}+V_x^{(p)}(-a(x-m_U^{(p)}(t))+m_{\phi}^{(p)}(t))-\frac{1}{4}(cV_x^{(p)}+\rho(x-m_U^{(p)}(t)))^2+\beta(x-m_U^{(p)}(t))^2\nonumber\\
		&\quad+\int_0^1(V^{(p)}(t,x-xz)-V^{(p)}(t,x))\nu(dz)
    \end{align}
    with the terminal condition $V^{(p)}(T,x)=\gamma(x-m_U^{(p)}(T))^2$ for all $x\in \mathbb{R}$.
   Consider an ansatz of the value function taking the form as:
    \begin{equation}\label{value}
		V^{(p)}(t,x)=A_p(t)(x-m_U^{(p)}(t))^2+B_p(t)(x-m_U^{(p)}(t))+C_p(t).
    \end{equation}
Plugging \eqref{value} into Eq.~\eqref{pde}, we have
    \begin{align}
		0&=
		\Bigg[A_p'(t)-c^2A_p^2(t)-\left(2a+\rho c-\int_0^1(z^2-2z)\nu(dz)\right)A_p(t)-\frac{1}{4}\rho^2+\beta\Bigg]x^2
		+\Bigg[B_p'(t)\nonumber\\&\quad-\Big[a+\frac{\rho c}{2}+c^2A_p(t)+\int_0^1z\nu(dz)\Big]B_p(t)-2A_p'(t)m_U^{(p)}(t)+2A_p(t)\Big[\left(2a+\rho c+c^2A_p(t)\right)m_U^{(p)}(t)\nonumber\\&\quad+m_{\phi}^{(p)}(t)\Big]+\left(\frac{\rho^2}{2}-2\beta +A_p(t)\right)m_U^{(p)}(t)-2A_p(t)(m_U^{(p)})'(t)\Bigg]x+\Bigg[C_p'(t)+A_p'(t)(m_U^{(p)})^2(t)\nonumber\\&\quad-B_p'(t)m_U^{(p)}(t)+\left(-2A_p(t)m_U^{(p)}(t)+B_p(t)\right)\left(a+\frac{\rho c}{2}-\int_0^1z\nu(dz)\right)m_U^{(p)}(t)-\frac{c^2}{4}\big(B_p(t)\nonumber\\&\quad-2A_p(t)m_U^{(p)}(t)\big)^2+(\beta-\frac{1}{4}\rho^2) (m_U^{(p)})^2(t)+2m_U^{(p)}(t)A_p(t)(m_U^{(p)})'(t)-B_p(t)(m_U^{(p)})'(t)\Bigg].
    \end{align}
    We thus obtain the following ODE system that $A_p(t)$, $B_p(t)$ and $C_p(t)$ should satisfy
	\begin{equation}\label{ODE}
		\left\{
		\begin{aligned}
			A_p'(t)&=c^2A_p^2(t)+(2a+\rho c-\int_0^1(z^2-2z)\nu(dz))A_p(t)+\frac{1}{4}\rho^2-\beta,\\
			B_p'(t)&=\left[a+c^2A_p(t)+\frac{\rho c}{2}+\int_0^1z\nu(dz)\right]B_p(t)-2A_p(t)m_U^{(p)}(t)\int_0^1(z^2-z)\nu(dz)\\
   &\quad+2A_p(t)(m_U^{(p)})'(t)-2A_p(t)m_{\phi}^{(p)}(t),\\
			C_p'(t)&= -\int_0^1z^2\nu(dz)A_p(t)(m_U^{(p)})^2(t)-2A_p(t)m_{\phi}^{(p)}(t)m_U^{(p)}(t)+\frac{1}{4}c^2B_p^2(t)\\
            &\quad+B_p(t)\left(m_U^{(p)}\right)'(t)
		\end{aligned}\right.
	\end{equation}
    with terminal conditions $A_p(T)=\gamma$ and  $B_p(T)=C_p(T)=0$.
    Under the condition $\rho^2< 4\beta$, we have $R>0$. Solving this Riccati equation to obtain that solution $t\mapsto (A_p(t),B_p(t),C_p(t))$ is given by \eqref{solution}. 
As a result, the  candidate optimal (feedback) control function and the solution of HJB equation have been finally obtained through the back substitutions of functions $A_p(t)$, $B_p(t)$ and $C_p(t)$ into \eqref{FOC} and \eqref{value}. We write $B_p(t, m_U^{(p)}, m_{\phi}^{(p)})$ and $C_p(t, m_U^{(p)}, m_{\phi}^{(p)})$ instead of $B_p(t)$ and $C_p(t)$ to highlight the dependence on $m_U^{(p)}$ and $m_{\phi}^{(p)}$.

We next prove that, the feedback strategy function $\theta^*(t,x)$ given by \eqref{optimalcontrol} is an optimal (feedback) control function. Note that $\delta^+>0$ and $\delta^-<0$. Then, $A_p(t)\geq 0$ for all $t\in [0, T]$, and hence $V_{xx}>0$ by using \eqref{value} and \eqref{solution}. For any $(t,x)\in[0,T]\times \mathbb{R}$, $s\in[t, T)$ and $\theta\in \mathcal{A}$, we have from It$\rm{\hat{o}}$'s formula that
	\begin{align*}
		V^{(p)}(s,U_s^{(p)})&=V^{(p)}(t,x)+\int_t^{s}V_t^{(p)}(v, U_v^{(p)})+V_x^{(p)}(v, U_v^{(p)})\left(a(m_U^{(p)}(v)-x)+m_{\phi}^{(p)}(v)+c\theta_v\right)dv\nonumber\\
  &\quad+\int_t^{s }\int_0^1(V^{(p)}(v, U_v^{(p)}-U_v^{(p)}z)-V^{(p)}(v, U_v^{(p)}))\nu(dz)dv.
	\end{align*}
Taking conditional expectation on both sides, we obtain that
\begin{align}\label{eq:ExVs}
\mathbb{E}_{t,x}[V^{(p)}(s, U_{s})]&=V^{(p)}(t,x)+\mathbb{E}_{t,x}\Bigg[\int_t^{s}V_t^{(p)}(v, U_v^{(p)})+V_x^{(p)}(v, U_v^{(p)})(a(m_U^{(p)}(v)-x)+m_{\phi}^{(p)}(v)+c\theta_v)dv\nonumber\\
&\quad+\int_t^{s}\int_0^1(V^{(p)}(v, U_v^{(p)}-U_v^{(p)}z)-V^{(p)}(v, U_v^{(p)}))\nu(dz)dv\Bigg]. 
\end{align}
Since the value function $V^{(p)}$ satisfies the HJB equation~\eqref{HJB}, we have, for any strategy  $\theta\in\mathcal{A}$,
\begin{align}\label{eq:NeqHJB}
&V_t^{(p)}(v, U_v^{(p)})+V_x^{(p)}(v, U_v^{(p)})(a(m_U^{(p)}(v)-x)+m_{\phi}^{(p)}(v)+c\theta_v \nonumber \\
&\quad +\int_0^1\big(V^{(p)}(v, U_v^{(p)}-U_v^{(p)}z)-V^{(p)}(v, U_v^{(p)})\big)\nu(dz)\nonumber\\
&\qquad\geq-\theta_v^2-\rho\theta_v\left(U_v^{(p)}-m_U^{(p)}(v)\right)-\beta\left|U_v^{(p)}-m_U^{(p)}(v)\right|^2.
\end{align}
Thus, plugging \eqref{eq:NeqHJB} into \eqref{eq:ExVs}, it is deduced that 
\begin{align*}
		\mathbb{E}_{t,x}\left[V^{(p)}(s, U_s^{(p)})\right]-V^{(p)}(t,x)\geq -\mathbb{E}_{t,x}\left[\int_t^s\theta_v^2+\rho\theta_v(U_v^{(p)}-m_U^{(p)}(v))+\beta\left|U_v^{(p)}-m_U^{(p)}(v)\right|^2dv\right].
\end{align*}
By sending $s$ to $T$ in the above display with the terminal condition $V^{(p)}(T,x)=\gamma(x-m_U^{(p)}(T))^2$, we arrive at, for any $\theta\in\mathcal{A}$,
\begin{equation*}
	V^{(p)}(t,x)\leq \mathbb{E}_{t,x}\left[\int_t^{T}\theta_s^2
 +\rho\theta_s(U_s^{(p)}-m_U^{(p)}(s))+\beta\left|U_s^{(p)}-m_U^{(p)}(s)\right|^2ds+\gamma\left|U_T-m_U^{(p)}(T)\right|^2\right],
	\end{equation*}
where the equality holds when one replaces strategy $\theta$ with $\theta^{*,(p)}$. This ends the proof of the lemma.
    \end{proof}
    
Lemma \ref{lem:bestresponse} provides the best response solution of the stochastic control problem  \eqref{dynamic}-\eqref{MFC} of the representative neuron. Under the optimum $\theta^{*,(p)}=(\theta_t^{*,(p)})_{t\in[0,T]}$ given in \eqref{optimalcontrol}, the corresponding membrane potential process  $U^{*,(p)}=(U_t^{*,(p)})_{t\in[0,T]}$ satisfies the dynamics given by
\begin{align}\label{u^*}
	U_t^{*,(p)}
	&= u-\int_0^t\left(a+c^2A_p(s)+\frac{\rho c}{2}\right)\left(U_s^{*,(p)}-m_U^{(p)}(s)\right)ds-\frac{c^2}{2}\int_0^tB_p\left(s, m_U^{(p)}, m_{\phi}^{(p)}\right)ds\nonumber\\
	&\quad-\int_0^t\int_0^1zU_{s-}^{*,(p)}N(ds, dz)+\int_0^tm_{\phi}^{(p)}(s)ds.
\end{align}

\subsection{The fixed point problem}
To characterize the mean field equilibrium as in Definition~\ref{def:MFE},  this subsection focuses on the fixed point problem arising from the consistency condition. The conditional consistency condition related to the MFG problem is given by, for any $p\in\mathcal{O}$,
	\begin{align}\label{consistency}
 m_{U}^{(p)}(t)=\mathbb{E}[U_t^{*,(p)}],\quad m_{\phi}^{(p)}(t)=\mathbb{E}\left[\int_0^1\phi(zU_t^{*,(p)})\nu(dz)\right],~~\forall t\in[0,T].
	\end{align}
Note that the connection strength function is of linear scenarios, i.e., $\phi(x)=kx+\ell$. In this scenario, once the first equation in \eqref{consistency} is true, the second holds automatically. In fact
    \begin{align*}
        \frac{1}{n}\sum_{j=1}^n\int_0^1\phi(zU_{t}^j)\nu(dz)=\frac{1}{n}\sum_{j=1}^n\int_0^1(kzU_{t}^j+\ell)\nu(dz)=k\left(\int_0^1z\nu(dz)\right)\left(\frac{1}{n}\sum_{j=1}^nU_{t}^j\right)+\ell.
    \end{align*}
Heuristically, $m_{\phi}^{(p)}(t)=k(\int_0^1z\nu(dz)) m_U^{(p)}(t)+\ell$ as $n \to \infty$. This yields that we can only consider the 1st consistency condition in \eqref{consistency}. It follows from \eqref{u^*} that the 1st condition in \eqref{consistency} can be read as follows: 
    \begin{align}\label{fix}
    	m_{U}^{(p)}(t)
    	=u-\int_0^tH(s, m_U^{(p)})ds+(k-1)\left(\int_0^1z\nu(dz)\right)\int_0^tm_U^{(p)}(s)ds+\ell t.
    \end{align}
Here, for $(t,m_U^{(p)})\in[0,T]\times\mathcal{C}_T^1$,     
\begin{align}
H(t, m_U^{(p)})&:=\frac{c^2}{2}B_p\left(t, m_U^{(p)}, m_{\phi}^{(p)}\right)\nonumber\\
&=c^2h_p(t)\left[-\gamma m_U^{(p)}(T)+\frac{A_p(t)}{h_p(t)}m_U^{(p)}(t)+\ell\int_t^T\frac{A_p(s)}{h_p(s)} ds \right.\\
     & \quad\left.+\int_t^T\frac{1}{h_p(s)}\left(\left(a+\frac{\rho c}{2}+k\int_0^1z\nu(dz)\right)A_p(s)+\left(\frac{\rho^2}{4}-\beta \right)\right)m_U^{(p)}(s)ds\right].\nonumber
    \end{align}
    with $h_p(t):=e^{-\int_t^T(a+c^2A_p(s)+\int_0^1z\nu(dz)+\frac{\rho c}{2})ds}$ for $t\in[0,T]$. Then, we have   
    
    \begin{proposition}\label{fixpoint}
 For any $p\in\mathcal{O}$, there exists a unique fixed point  $m_U^{*,(p)}= (m_U^{*,(p)}(t))_{t\in[0,T]} \in \mathcal{C}^1_T$ to the fixed point problem \eqref{fix}. Moreover, the conditional mean field equilibrium is given by  \begin{equation}\label{mfe}
    \theta_t^{*,(p)}=\left(-cA_p(t)-\frac{\rho}{2}\right)(U_t^{*,(p)}-m_U^{*,(p)}(t))-\frac{c}{2}B_p\left(t,m_U^{*,(p)}, m_{\phi}^{*,(p)}\right),
    \end{equation}
    where $U^{*,(p)}=(U_t^{*,(p)})_{t\in[0,T]}$ is the membrane potential process under optimal control $\theta^{*,(p)}$, and $m_{\phi}^{*,(p)}(t)=k\left(\int_0^1z\nu(dz)\right) m_U^{*,(p)}(t)+\ell$.
    \end{proposition}
    
 Nonetheless, deriving the existence of a fixed point directly on $\mathcal{C}^1_T$ proves to be not a direct task. Consequently, we divide the proof into two distinct stages. Initially, we aim to identify a fixed point $m_U^{(p)} \in \mathcal{C}_T:= C([0,T];\mathbb{R})$ for \eqref{fix}. Subsequently, we demonstrate that this solution denoted by $m_U^{*,(p)}=(m_U^{*,(p)}(t))_{t\in[0,T]}$ indeed resides within $\mathcal{C}^1_T$. To do it, we first have the following auxiliary lemma:
    \begin{lemma}\label{lem:C-fixed}
For any $p\in\mathcal{O}$, there exists a unique fixed point $m_U^{*,(p)}= (m_U^{*,(p)}(t))_{t\in[0,T]} \in \mathcal{C}_T$ to \eqref{fix}.
    \end{lemma}
    
    \begin{proof}
Let us first define the following mapping $\Phi$ on $\mathcal{C}_T$ by, for $(t, m_U)\in[0,T]\times \mathcal{C}_T$,
\begin{align}\label{map1}
\Phi(t, m_U)&:= u-
    \int_0^tH(s, m_U)ds+(k-1)\left(\int_0^1z\nu(dz)\right)\int_0^tm_U(s)ds+\ell t.
\end{align}
It is not difficult to observe from \eqref{map1} that $\Phi(\cdot,m_U)\in\mathcal{C}_T^1\subset\mathcal{C}_T$ for all $m_U\in\mathcal{C}_T$. Moreover, we have, for any $m_U^1, m_U^2\in \mathcal{C}_T$,
\begin{align*}
&\left| H(t, m_U^1)-H(t, m_U^2)\right|=\bigg|-c^2h_p(t)\gamma(m_U^1(T)-m_U^2(T))+c^2A_p(t)(m_U^1(t)-m_U^2(t))\nonumber\\
    &\quad +c^2\int_t^T\frac{h_p(t)}{h_p(s)}\left(A_p(s)\left(a+\frac{\rho c}{2}+k\int_0^1z\nu(dz)\right)+\left(\frac{\rho^2}{4}-\beta \right)\right)(m_U^1(s)-m_U^2(s))ds\bigg|\nonumber\\
    &~~\leq M_1(T)\left\| m_U^1-m_U^2\right\|_T,
    \end{align*}
where $\|f\|_T:=\sup_{t\in[0,T]}|f(t)|$ for any $f\in\mathcal{C}_T$, and the mapping 
\begin{align}\label{eq:M1t}
M_1(t)&:=\gamma c^2\left\|h_p\right\|_T+c^2\left\|A_p\right\|_T+\left\|h_p\right\|_Tc^2\left\|\frac{A_p\left(a+\frac{\rho c}{2}+k\int_0^1z\nu(dz)\right)+\left(\frac{\rho^2}{4}-\beta \right)}{h_p}\right\|_Tt.  
\end{align}
Then, it follows from \eqref{map1} that
\begin{align*}
\left|\Phi(t, m_U^1)-\Phi(t, m_U^2)\right|
&\leq \left[M_1(T)t+\left|(k-1)\int_0^1z\nu(dz)\right| t\right]\left\|m_U^1-m_U^2\right\|_T.
\end{align*}
Introduce $M(T):=M_1(T)T+\left|(k-1)\left(\int_0^1z\nu(dz)\right)\right| T$. Then, $\|\Phi(\cdot, m_U^1)-\Phi(\cdot, m_U^2)\|_T\leq M(T)\|m_U^1-m_U^2\|_T$ for all $m_U^1, m_U^2\in \mathcal{C}_T$. Note that $T \mapsto M(T)$ is continuous on $\mathbb{R}$ and satisfies $\lim_{T\to0} M(T) = 0$. Thus, we may choose a sufficiently small $t_1\in(0,T]$ such that $M(t_1) \in(0, 1)$. Hence, there exists a unique fixed point of $\Phi$ on $\mathcal{C}_T$ with $t\in[t_0, t_1]$ and $t_0 = 0$, since $\Phi$ is a contraction mapping on $\mathcal{C}_{t_1}$.  Similarly, we can examine the existence of a unique fixed point of $\Phi$ with  $t \in[t_1, t_2]$ for some $t_2 > t_1$ small enough. By going through this process again, we can conclude there exists a unique fixed point $m_U^{*,(p)}$ of $\Phi$ on $[0, T]$.
\end{proof}

We are now at position to prove Proposition \ref{fixpoint}.
    
\begin{proof}[Proof of Proposition \ref{fixpoint}.]
    We aim to apply the implicit function theorem (c.f. Theorem 15.1 and Corollary 15.1 of \cite{D13}) to show that the fixed point $m_U^{*,(p)}\in\mathcal{C}_T$ in Lemma \ref{lem:C-fixed} is also in $\mathcal{C}^1_T$. For this purpose, define the following mapping $F: [0,T]\times \mathcal{C}_T\to \mathbb{R}$ as:
    \begin{align}\label{eq:FtMu}
     F(t, m_U):=\Phi(t, m_U)-m_U(t),\quad \forall (t, m_U)\in[0,T]\times\mathcal{C}_T,   
    \end{align}
where $\Phi(t, m_U)$ is defined by \eqref{map1}. Obviously $F(t, m_U^{*,(p)})=0$. Firstly, the mapping $F$ is continuous on $[0,T]\times\mathcal{C}_T$. In fact, for any $m_U^1, m_U^2\in \mathcal{C}_T$ and $s, t\in[0,T]$, we have
    \begin{align*}
    \left|F(s, m_U^1)-F(t, m_U^2)\right|
    &\leq \left|\Phi(s, m_U^1)-\Phi(s, m_U^2)\right|+\left|\Phi(s, m_U^2)-\Phi(t, m_U^2)\right|+\left|m_U^1(s)-m_U^2(s)\right|\\
    &\quad+\left|m_U^2(s)-m_U^2(t)\right|\\
    &\leq M(T)\left\|m_U^1-m_U^2\right\|_T+M_2(T)|s-t|+\left\|m_U^1-m_U^2\right\|_T+\left|m_U^2(s)-m_U^2(t)\right|,
    \end{align*} 
    where $M_2(T)=M_1(T)\|m_U^2\|_T+\left|(k-1)\int_0^1z\nu(dz)\right| \|m_U^2\|_T+\left\|\frac{A_p}{h_p}\right\|_Tc^2\ell T+\ell$.

We next prove that $F$ is continuously differentiable on $[0,T]\times\mathcal{C}_T$. To do it, denote $F_t(t,m_U)$ and $F_{m_U}(t,m_U)$ as the (Frech\'{e}t) derivatives of mapping $F$ w.r.t $t$ and $m_U$, respectively. Additionally,  $\Phi_t(t, m_U)$ and $\Phi_{m_U}(t, m_U)$ represent the (Frech\'{e}t) derivatives of mapping $\Phi$ w.r.t $t$ and $m_U$. Then, it holds that
\begin{align}\label{Frecht}
 F_t(t, m_U)=\Phi_t(t,m_U),\quad F_{m_U}(t,m_U)=\Phi_{m_U}(t,m_U)-I, 
\end{align} 
where $I$ is the identity operator, $\Phi_t(t,m_U)=-H(t, m_U)+(k-1)(\int_0^1z\nu(dz))m_U(t)+\ell$ and $\Phi_{m_U}(t,m_U) x =-\int_0^tH(s, x)ds-\ell c^2\int_0^th_p(s)\int_s^T\frac{A_p(v)}{h_p(v)}dvds+(k-1)\int_0^1z\nu(dz)\int_0^{t} x(s)ds$ for $x\in \mathcal{C}_T$. Consequently, for any $m_U^1, m_U^2\in \mathcal{C}_T$ and $s, t\in[0,T]$, we have
     \begin{align*}
    &|F_t(s, m_U^1)-F_t(t, m_U^2)|\leq |F_t(s, m_U^1)-F(t, m_U^1)|+|F_t(t, m_U^1)-F(t, m_U^2)|\\
    &\leq \left|H(t, m_U^2)-H(t, m_U^1)\right|+\left|H(t,m_U^1)-H(s, m_U^1)\right|+\left|(k-1)\int_0^1z\nu(dz)\right|\left(\left|m_U^1(s)-m_U^2(s)\right|\right.\\&\quad\left.+\left|m_U^2(s)-m_U^2(t)\right|\right)
    \\&\leq M_1(T)\left\|m_U^1-m_U^2\right\|_T+\left|H(t,m_U^1)-H(s, m_U^1)\right|+\left|(k-1)\int_0^1 z\nu(dz) \right|\left(\left\|m_U^1-m_U^2\right\|_T\right.\\&\quad\left.+\left|m_U^2(s)-m_U^2(t)\right|\right),
    \end{align*}
   where the following estimate hods true:
    \begin{align*}
    &|H(t, m_U^1)-H(s, m_U^1)|\leq \gamma m_U^1(T)c^2|h_p(s)-h_p(t)|+c^2|A_p(t)-A_p(s)|\|m_U^1\|_T+c^2\|A_p\|_T\\&\quad\times|m_U^1(t)-m_U^1(s)|+\bigg\|c^2 \ell\frac{A_p}{h_p}+c^2\frac{1}{h_p}\left(\left(a+\frac{\rho c}{2}+k\int_0^1z\nu(dz)\right)A_p+\left(\frac{\rho^2}{4}-\beta \right)\right) m_U^1\bigg\|_T T\\&\quad \times|h_p(t)-h_p(s)|+\bigg\|c^2 \ell\frac{A_p}{h_p}+c^2\frac{1}{h_p}\left(\left(a+\frac{\rho c}{2}+k\int_0^1z\nu(dz)\right)A_p+\left(\frac{\rho^2}{4}-\beta \right)\right) m_U^1\bigg\|_T\|h_p\|_T\\&\quad \times|s-t|
    \\&\leq C_T\left\{|h_p(s)-h_p(t)|+|A_p(t)-A_p(s)|+|m_U^1(t)-m_U^1(s)|+|s-t|\right\}.
    \end{align*} 
Here, $C_T>0$ is a constant depending on $T$ only. 
Then, it follows from the continuity of $t\mapsto h_p(t)$ and $t\mapsto A_p(t)$ that $F_t(t, m_U)$ is continuous on $[0,T]\times \mathcal{C}_T$. Note that, for $s, t\in[0, T]$ and $m_U^1, m_U^2 \in \mathcal{C}_T$,
\begin{align*}
\left|F_{m_U}(s,m_U^1)x-F_{m_U}(t,m_U^2)x\right|&\leq \left[M_1(T)+\Big|(k-1)\int_0^1z\nu(dz)\Big|\right]\|x\|_T\times|s-t|+|x(s)-x(t)|.
\end{align*}
This yields that $F_{m_U}(t,m_U)$ is continuous on $[0,T]\times \mathcal{C}_T$. 
For fixed $(t, m_U)\in [0,T]\times \mathcal{C}_T$, $F_{m_U}(t,m_U)$ is a mapping from $\mathcal{C}_{T}$ to $\mathbb{R}$. 
Then, lastly, we need to prove that for any $t\in [0, T]$, $F_{m_U}^{-1}(t,m_U^{*,(p)})\in \mathcal{L}(\mathbb{R}; \mathcal{C}_T)$, where $\mathcal{L}(\mathbb{R}; \mathcal{C}_T)$ is the set of bounded linear operators from $\mathbb{R}$ to $\mathcal{C}_T$. Here, we recall that the fixed point $m_U^{*,(p)}\in\mathcal{C}_T$ is given in Lemma \ref{lem:C-fixed}. Then, it suffices to show that for any $t\in [0, T]$, $F_{m_U}(t,m_U^{*,(p)})$ is one-to-one and onto $\mathbb{R}$. We will show that, for any $\alpha \in \mathbb{R}$ and $t\in [0, T]$, there exists an unique solution $x \in \mathcal{C}_T$ to the equation $F_{m_U}(t, m_U^{*,(p)})x=\alpha$, or equivalently $x=\Phi_{m_U}(t,m_U^{*,(p)})x -\alpha$. 
Note that, for the mapping $\Phi_{m_U}(\cdot,m_U^{*,(p)}):\mathcal{C}_T\to \mathcal{C}_T$, we have
\begin{align*}
\sup_{t\in[0,T]}\left|\Phi_{m_U}(t,m_U^{*,(p)})x_1-\alpha-\big(\Phi_{m_U}(t,m_U^{*,(p)})x_2-\alpha\big)\right|\leq M(T) \Vert x_1-x_2\Vert_T,\quad\forall x_1,x_2\in\mathcal{C}_T,   
\end{align*}
where $T\mapsto M(T)$ is given in the proof of Lemma~\ref{lem:C-fixed}, which satisfies that $T\mapsto M(T)$ is continuous and $\lim_{T\to 0}M(T)=0$. Following the similar proof of Lemma \ref{lem:C-fixed}, we may choose a sufficiently small $t_1\in(0,T]$ such that $M(t_1) \in(0, 1)$. Hence, $\Phi_{m_U}(\cdot,m_U^{*,(p)})$ is a contraction mapping on $\mathcal{C}_{t_1}$.  By going through this process again, we can conclude there exists a unique solution to $x=\Phi_{m_U}(t,m_U^{*,(p)})x-\alpha$ for any $t\in [0, T]$, which yields that $F_{m_U}(t,m_U^{*,(p)})$ maps onto $\mathbb{R}$. 
In the case $\alpha=0$, we can obtain that there exits a unique solution to $F_{m_U}(t, m_U^{*,(p)})x=0$ in $\mathcal{C}_T$. Since $x = 0$ is a solution to $F_{m_U}(t,m_U^{*,(p)})x=0$ for any $t\in[0,T]$. Hence,  we have $x_0\equiv 0$, which implies that $F_{m_U}(t,m_U^{*,(p)})$ is one-to-one. By virtue of the open mapping theorem (c.f. Proposition 7.8 of \cite{D13}), we have $(F_{m_U})^{-1}(t, m_U^{*,(p)}) \in \mathcal{L}(\mathcal{C}_T;\mathbb{R})$. By applying the implicit function theorem, for every $t\in [0,T]$, $m_U^{*,p}$ is continuously differentiable on $[0, T]$, thus $m_U^{*,p}\in \mathcal{C}^1_T$. According to Definition \ref{def:MFE}, we can claim that \eqref{mfe} is the conditional mean field equilibrium, which completes the proof.
    \end{proof}

\section{Approximate Nash Equilibrium}\label{sec:approximateNE}

In this section, we aim to construct the approximate Nash equilibrium for the finite number of neurons in the coupling control system \eqref{ndynamics}-\eqref{objective}. 

Let $\boldsymbol{p}=(\boldsymbol{u},\boldsymbol{a}, \boldsymbol{c})\sim\mu$ be a random type vector independent of $N(dt,dz)$ under the filtered probability space $(\Omega, \mathcal{F}, \mathbb{F}, \mathbb{P})$. We define that, for $t\in[0,T]$,
\begin{align}\label{eq:mstar}
m_U^*(t):=\mathbb{E}_{\mu}\left[m_U^{*,(\boldsymbol{p})}(t)\right]=\int_{\mathcal{O}}m_U^{*,(p)}(t)\mu(dp),\quad m_{\phi}^*(t):=k\left(\int_0^1z\nu(dz)\right)m_U^*(t)+\ell,    
\end{align}
where for $p\in\mathcal{O}$, $m_U^{*,(p)}=(m_U^{*,(p)}(t))_{t\in[0,T]}\in\mathcal{C}_T^1$ is the fixed point provided in Proposition \ref{fixpoint}.
By virtue of Proposition~\ref{fixpoint} for our MFG problem,  we construct the following feedback control strategy of neuron $i$, for $i=1, \ldots, n$,
    \begin{align}\label{ntheta}
    \theta_t^{*,i}=\theta^{*,i}(t, U_t^{*,i})&:=\left(-c_iA_i(t)-\frac{\rho}{2}\right)(U_t^{*,i}-m_U^*(t))-\frac{c_i}{2}B_i\left(t, m_U^*, m_{\phi}^*\right),
    \end{align}
where  $m_U^*,m_{\phi}^*$ are defined by \eqref{eq:mstar}, and the coefficients are given by, for $m_U,m_{\phi}\in\mathcal{C}_T^1$,
    \begin{equation}\label{npara1}	
    \left\{
    \begin{aligned}
    &A_i(t)=\frac{(\frac{\rho^2}{4}-\beta)(e^{(\delta_i^+-\delta_i^-)(T-t)}-1)-\gamma( \delta_i^+e^{(\delta_i^+-\delta_i^-)(T-t)})-\delta_i^-)}{\delta_i^-e^{(\delta_i^+-\delta_i^-)(T-t)}-\delta_i^+-c_i^2\gamma(e^{(\delta_i^+-\delta_i^-)(T-t)}-1)},\\
	&B_i(t, m_U, m_{\phi}):=2h_i(t)\Bigg[-\int_t^T\frac{A_i(s)}{h_i(s)}\left(m_U'(s)-m_{\phi}(s)-m_U(s)\int_0^1(z^2-z)\nu(dz)\right)ds\Bigg],\\[0.4em]
	&h_i(t):=e^{-\int_t^T(a_i+c_i^2A_i(s)+\int_0^1z\nu(dz)+\frac{\rho c_i}{2})ds}
    \end{aligned}
    \right.
    \end{equation}
    with $R_i:=(a_i+\frac{\rho c_i}{2}-\frac{1}{2}\int_0^1(z^2-2z)\nu(dz))^2-c_i^2(\frac{1}{4}\rho^2-\beta)>0$ and $\delta_i^{\pm}:=-(a_i+\frac{\rho c_i}{2}-\frac{1}{2}\int_0^1(z^2-2z)\nu(dz))\pm\sqrt{R_i}$. Here, $U^{*,i}= (U_t^{*,i})_{t\in[0,T]}$ is the membrane potential process under the control strategy $\theta^{*,i}$ given by \eqref{ntheta}, and hence it satisfies the dynamics given by
    \begin{align}\label{ui*}
    U_t^{*,i}&=u_i-a_i\int_0^t(U_s^{*,i}-\bar{U}_s^*)ds-\int_0^t\left(c_i^2A_i(s)+\frac{\rho c_i}{2}\right)(U_s^{*,i}-m_U^*(s))-\frac{c_i^2}{2}\int_0^tB_i(s,m_U^*,m_{\phi}^*)ds\nonumber\\&\quad-\int_0^t\int_0^1zU_{s-}^{*,i}N^i(ds, dz)
    +\frac{1}{n}\sum_{j\neq i}\int_0^t\int_0^1\phi(zU_{s-}^{*,j})N^j(ds,dz).
    \end{align}
    Here, $\bar{U}_t^*:=\frac{1}{n}\sum_{i=1}^nU_t^{*, i}$ denotes the average of membrane potential of neuron system at time $t\in[0,T]$.

    Next, we introduce the main result of this section.
    \begin{theorem}\label{theorem}
    Let Assumption \ref{Ap} hold. Then the admissible control vector  $\boldsymbol{\theta}^{*}=(\theta^{*, 1}_t,\ldots, \theta^{*, n}_t)_{t\in[0,T]}$ given by \eqref{ntheta} is an $\epsilon_n$-Nash equilibrium to the $n$-player game problem, namely 
    \begin{equation}\label{ANE}
    J_i(\boldsymbol{\theta}^*)-\epsilon_n\leq \inf_{\theta^i\in\mathcal{A}}J_i(\theta^i, \boldsymbol{\theta}^{*, -i}) , \quad \forall i=1, \ldots, n
    \end{equation}
   with $\epsilon_n>0$ satisfying $\lim_{n\to\infty}\epsilon_n=0$.
    \end{theorem}
   
To prove Theorem \ref{theorem}, we first present the following auxiliary results, the proof of which is in the appendix.
    \begin{lemma}\label{lemma2}
    Let Assumption \ref{Ap} hold. Then, it holds that
    \begin{itemize}
    
    \item[{\rm(i)}] for $i=1,\ldots,n$, and any strategy $\theta^i\in \mathcal{A}$ satisfying $\|\theta^i\|_{2,T}\leq\sqrt{\kappa}$ for some positive $\kappa$ depending on $T$ only, there is $C>0$ independent of $(i,n)$ s.t. $\mathbb{E}[\sup_{t\in [0,T]}\vert U_t^i\vert^2]\leq C$.
    
    \item[{\rm(ii)}] for $i=1,\ldots, n$, there is $C>0$ independent of $(i,n)$ s.t. $\mathbb{E}[\sup_{t\in[0,T]}\vert U_t^{*,i}\vert^2]\leq C$.
   
    \item[{\rm(iii)}] for $m_U^*$ given in \eqref{eq:mstar}, $\lim_{n\to \infty}\mathbb{E}[\vert\bar{U}_t^*-m_U^*(t)\vert^2]=0$ for any $t\in[0,T]$.
    \end{itemize}
    \end{lemma}

Based on the preparations provided above, we are ready to prove Theorem \ref{theorem}.

    \begin{proof}[Proof of Theorem \ref{theorem}]
    	
Let us first introduce that, for $i=1,\ldots,n$,         
    \begin{equation}\label{checkU}
    		\left\{\begin{aligned}
    			\check{U}_t^i&=u_i-a_i\int_0^t\left(\check{U}_s^i-\frac{1}{n}\sum_{l=1}^n\check{U}_s^l\right)ds+\int_0^tc_i\theta_s^ids-\int_0^t\int_0^1z\check{U}_{s-}^iN^i(ds, dz)\\&\quad+\frac{1}{n}\sum_{l\neq i}\int_0^t\int_0^1\phi(z\check{U}_s^l)N^l(ds, dz),\\
    			\check{U}_t^j&=u_j-a_j\int_0^t\left(\check{U}_s^j-\frac{1}{n}\sum_{l=1}^n\check{U}_s^l\right)ds+\int_0^tc_j\check{\theta}_s^{*,j}ds-\int_0^t\int_0^1z\check{U}_{s-}^jN^j(ds,dz)\\&\quad+\frac{1}{n}\sum_{l\neq j}\int_0^t\int_0^1\phi(z\check{U}_s^l)N^l(ds,dz), \qquad j\neq i,
    		\end{aligned}
    		\right.
    	\end{equation}
        where $\check{\theta}^{*,j}_t=\theta^{*,j}(t, \check{U}_t^j)=\left(-c_jA_j(t)-\frac{\rho}{2}\right)(\check{U}_t^{*,j}-m_U^*(t))-\frac{c_j}{2}B_j(t, m_U^*, m_{\phi}^*)$ for $j\neq i$, while $\theta^i\in\mathcal{A}$. Following a similar argument as in Lemma \ref{lemma2}-({\romannumeral1}), there exists a constant $C>0$ depending on $T>0$ only such that, for all $t\in[0,T]$,
        \begin{equation}
        	\mathbb{E}\left[\vert \check{U}_t^i\vert^2\right]\leq C\left\{1+\left\|\theta^i\right\|_{2,T}^2\right\}, \quad
        	\mathbb{E}\left[\vert \check{U}_t^j\vert^2\right]\leq C\left\{1+\frac{1}{n}\left\|\theta^i\right\|_{2,T}^2\right\}.
        \end{equation} 
Then, we have that
\begin{align*}
J_i(\theta^i, \boldsymbol{\theta}^{*,-i}) &=\mathbb{E}\left[\int_0^T\left((\theta_t^i)^2+\rho\theta_t^i\left(\check{U}_t^i-\frac{1}{n}\sum_{j=1}^n\check{U}_t^j\right)+\beta\left\vert\check{U}_t^i-\frac{1}{n}\sum_{j=1}^n\check{U}_t^j\right\vert^2\right)dt \right.\\
&\quad \left.
+\gamma\left\vert\check{U}_T^i-\frac{1}{n}\sum_{j=1}^n \check{U}_T^j\right\vert^2\right].
\end{align*}
Note that, it holds that 
\begin{align*}
    &J_i(\theta^i, \boldsymbol{\theta}^{*,-i})\\&=\mathbb{E}\left[\int_0^T\left((\theta_t^i)^2+\rho\theta_t^i\left(\check{U}_t^i-\frac{1}{n}\sum_{j=1}^n\check{U}_t^j\right)+\beta\left\vert\check{U}_t^i-\frac{1}{n}\sum_{j=1}^n\check{U}_t^j\right\vert^2\right)dt +\gamma\left\vert\check{U}_T^i-\frac{1}{n}\sum_{j=1}^n\check{U}_T^j\right\vert^2\right]\\
    &\geq \mathbb{E}\left[\int_0^T\left((\theta_t^i)^2-\rho\vert\theta_t^i\vert\left\vert \check{U}_t^i-\frac{1}{n}\sum_{j=1}^n\check{U}_t^j\right\vert+\beta\left\vert\check{U}_t^i-\frac{1}{n}\sum_{j=1}^n\check{U}_t^j\right\vert^2\right)dt\right]\\
    &\geq \mathbb{E}\left[\int_0^T\left(\beta\left(\left\vert \check{U}_t^i-\frac{1}{n}\sum_{j=1}^n\check{U}_t^j\right\vert-\frac{\rho}{2\beta}\vert\theta_t^i\vert\right)^2+\left(1-\frac{\rho^2}{4\beta}\right)\vert\theta_t^i\vert^2\right)dt\right]\\
    &\geq\left(1-\frac{\rho^2}{4\beta}\right)\left\|\theta^i\right\|_{2,T}^2.
\end{align*}
This yields that $J_i(\theta^i, \boldsymbol{\theta}^{*,-i})$ is coercive in the sense that $J_i(\theta^i, \boldsymbol{\theta}^{*,-i})\to \infty $ as $\|\theta^i\|_{2,T}\to\infty$. Hence, it is enough to show $J_i(\theta^i, \boldsymbol{\theta}^{*,-i})\leq J(\boldsymbol{\theta}^*)$ under the condition that $\theta^i\in\mathcal{A}$ satisfying 
\begin{align}\label{eq:boundednesstheta2}
\left\|\theta^i\right\|_{2,T}^2=\mathbb{E}\left[\int_0^T\vert\theta_t^i\vert^2dt\right]\leq \kappa,~~\text{for}~\kappa>0.    
\end{align} 
In what follows, we introduce the following optimal control problem given by
    	\begin{align*}
    		\inf_{\theta^i\in\mathcal{A}}&\bar{J}_i(\theta^i; m_U^*):=\mathbb{E}\left[\int_0^T\left((\theta_t^i)^2+\rho\theta_t^i(\check{U}_t^i-m_U^*(t))+\beta\vert\check{U}_t^i-m_U^*(t)\vert^2\right)dt+\gamma\left\vert\check{U}^i_T-m_U^*(T) \right\vert^2\right],
    	\end{align*}
where $\check{U}=(\check{U}_t^i)_{t\in[0,T]}$ is defined in \eqref{checkU}, and $(m_U^*(t))_{t\in[0,T]}$ is given in \eqref{eq:mstar}. Then, it holds that
    	\begin{align}\label{A1}
    		&\inf_{\theta^i\in\mathcal{A}^i}J_i(\theta^i, \boldsymbol{\theta}^{*, -i})-J_i(\boldsymbol{\theta}^*)
    		\\
&=\left(\inf_{\theta^i\in\mathcal{A}^i}J_i(\theta^i, \boldsymbol{\theta}^{*,-i})-\inf_{\theta^i\in\mathcal{A}^i}\bar{J}_i(\theta^i; m_U^*) \right)+\inf_{\theta^i\in\mathcal{A}^i}\bar{J}_i(\theta^i; m_U^*)-\bar{J}_i(\check{\theta}^{*,i}; m_U^*)+\bar{J}_i(\check{\theta}^{*,i}; m_U^*)-J_i(\boldsymbol{\theta}^*)\nonumber\\
    		&\geq\inf_{\theta^i\in\mathcal{A}^i}\left(J_i(\theta^i, \boldsymbol{\theta}^{*, -i})-\bar{J}_i(\theta^i; m_U^*)\right)+\left( \inf_{\theta^i\in\mathcal{A}^i}\bar{J}_i(\theta^i; m_U^*)-\bar{J}_i(\check{\theta}^{*,i}; m_U^*)\right)+\left(\bar{J}_i(\check{\theta}^{*,i}; m_U^*)-J_i(\boldsymbol{\theta}^*)\right).\nonumber
    	\end{align}
First of all, we count on the first term on the RHS of \eqref{A1} and have that
    	\begin{align*}
    		&\left\vert J_i(\theta^i, \boldsymbol{\theta}^{*, -i})-\bar{J}_i(\theta^i; m_U^*)\right\vert \nonumber\\&\leq \rho\mathbb{E}\left[\int_0^T\left\vert\theta_t^i\left(m_U^*(t)-\frac{1}{n}\sum_{j=1}^n\check{U}_t^j\right)\right\vert dt\right]+\beta\mathbb{E}\left[\int_0^T\left\vert\left(\frac{1}{n}\sum_{j=1}^n\check{U}_t^j\right)^2-(m_U^*(t))^2\right\vert dt\right]\nonumber\\&\quad+2\beta\mathbb{E}\left[\int_0^T\left\vert \check{U}_t^i\left(m_U^*(t)-\frac{1}{n}\sum_{j=1}^n\check{U}_t^j\right)\right\vert dt\right]+\gamma\mathbb{E}\left[\left\vert\left(\frac{1}{n}\sum_{j=1}^n\check{U}_T^j\right)^2-(m_U^*(T))^2\right\vert \right]\nonumber\\&\quad+2\gamma \mathbb{E}\left[\left\vert \check{U}_T^i\left(m_U^*(T)-\frac{1}{n}\sum_{j=1}^n\check{U}_T^j\right)\right\vert \right]=: E_1^i+E_2^i+E_3^i+E_4^i+E_5^i.
    	\end{align*}
For the term $E_1^i$, by using Cauchy-Schwartz's inequality and \eqref{eq:boundednesstheta2}, we have
    	\begin{align*}
E_1^i&=\rho\mathbb{E}\left[\int_0^T\left\vert\theta_t^i\left(m_U^*(t)-\frac{1}{n}\sum_{j=1}^n\check{U}_t^j\right)\right\vert dt\right]\leq \rho\left\{\left\|\theta^{i}\right\|_{2,T}^2\int_0^T\mathbb{E}\left[\left\vert \frac{1}{n}\sum_{j=1}^n\check{U}_t^j-m_U^*(t)\right\vert ^2 \right]dt\right\}^{\frac{1}{2}}\\
    		&\leq C\left\{\int_0^T\mathbb{E}\left[\left\vert \frac{1}{n}\sum_{j=1}^n\check{U}_t^j-m_U^*(t)\right\vert ^2 \right]dt\right\}^{\frac{1}{2}},
    	\end{align*}
where 
    	\begin{align*}
    		\mathbb{E}\left[\left\vert\frac{1}{n}\sum_{j=1}^n\check{U}_t^j-m_U^*(t)\right\vert^2 \right]&\leq 2\mathbb{E}\left[\left\vert\frac{1}{n}\sum_{j=1}^n\check{U}_t^j-\bar{U}_t^*\right\vert^2 \right]+2\mathbb{E}\left[\left\vert\bar{U}_t^*-m_U^*(t)\right\vert^2 \right].
    	\end{align*}
On the other hand,  by recalling \eqref{ui*} and \eqref{checkU}, we have
    	\begin{equation*}
    		\left\{
    		\begin{aligned}
    		\check{U}_t^i-U_t^{*,i}&=a_i\int_0^t(U_s^{*,i}-\check{U}_s^i)ds+\frac{1}{n}\sum_{l=1}^na_i\int_0^t(\check{U}_s^l-U_s^{*,l})ds+c_i\int_0^t(\theta_s^i-\theta_s^{*,i})ds\\&\quad-\int_0^t\int_0^1(\check{U}_s^i-U_s^{*,i})N^i(ds,dz)+\frac{1}{n}\sum_{l\neq i}\int_0^t\int_0^1kz(\check{U}_s^l-U_s^{*,l})N^l(ds,dz),\\
    		\check{U}_t^j-U_t^{*,j}&=a_j\int_0^t(U_s^{*,j}-\check{U}_s^j)ds+\frac{1}{n}\sum_{l=1}^na_j\int_0^t(\check{U}_s^l-U_s^{*,l})ds+c_j\int_0^t(\check{\theta}_s^{*,j}-\theta_s^{*,j})ds\\&\quad-\int_0^t\int_0^1(\check{U}_s^j-U_s^{*,j})N^j(ds,dz)+\frac{1}{n}\sum_{l\neq j}\int_0^t\int_0^1kz(\check{U}_s^l-U_s^{*,l})N^l(ds,dz).\end{aligned}\right.
    	\end{equation*}
By using  a straightforward calculation as in Lemma \ref{lemma2}, it holds that
        \begin{equation*}
        	\left\{\begin{aligned}
        	\mathbb{E}\left[\vert\check{U}_t^i-U_t^{*,i}\vert^2\right]&\leq C\left\{\int_0^t\mathbb{E}\left[\vert\check{U}_s^i-U_s^{*,i}\vert^2\right]ds+\left\|\theta^i-\theta^{*,i}\right\|_{2,t}^2+\int_0^t\frac{1}{n}\sum_{l=1}^n\mathbb{E}\left[\vert\check{U}_s^l-U_s^{*,l}\vert^2\right]ds\right\},\\
        	\mathbb{E}\left[\vert\check{U}_t^j-U_t^{*,j}\vert^2\right]&\leq C\left\{\int_0^t\mathbb{E}\left[\vert\check{U}_s^j-U_s^{*,j}\vert^2\right]ds+\int_0^t\frac{1}{n}\sum_{l=1}^n\mathbb{E}\left[\vert\check{U}_s^l-U_s^{*,l}\vert^2\right]ds\right\}, \quad  j\neq i.	
        	\end{aligned}\right.
        \end{equation*}
As a result, by averaging the above equality form  $l=1,\ldots,n$, we arrive at,
    	\begin{align}\label{eq:temp-U}
\frac{1}{n}\sum_{l=1}^n\mathbb{E}\left[\vert\check{U}_t^l-U_t^{*,l}\vert^2\right]\leq C\left\{\frac{1}{n}\left\|\theta^i-\theta^{*,i}\right\|_{2,t}^2+\int_0^t\frac{1}{n}\sum_{l=1}^n\mathbb{E}\left[\vert\check{U}_s^l-U_s^{*,l}\vert^2\right]ds\right\}.
    	\end{align} 
Note that $\|\theta^i-\theta^{*,i}\|_{2,t}^2\leq 2(\|\theta^{i}\|_{2,T}^2+\|\theta^{*,i}\|_{2,T}^2)$. Then, Lemma \ref{lemma2} and \eqref{eq:boundednesstheta2} imply that $\frac{1}{n}\|\theta^i-\theta^{*,i}\|_{2,t}^2\rightarrow 0$ as $n\rightarrow \infty$. Applying the Gronwall's lemma to \eqref{eq:temp-U} and combining Jensen inequality, it deduces that, as $n\to\infty$, 
    	\begin{align*}
    		\mathbb{E}\left[\left\vert\frac{1}{n}\sum_{j=1}^n\check{U}_t^j-\bar{U}_t^*\right\vert^2\right]&=\mathbb{E}\left[\left\vert\frac{1}{n}\sum_{j=1}^n(\check{U}_t^j-U_t^{*,j})\right\vert^2\right]\leq \frac{1}{n}\sum_{j=1}^n\mathbb{E}\left[\vert\check{U}_t^j-U_t^{*,j}\vert^2\right]\to 0.
    	\end{align*}
It follows from \eqref{eq:boundednesstheta2} that $\mathbb{E}[\sup_{t\in [0,T]}\vert \check{U}_t^i\vert^2]\leq C$ for some finite $C>0$ depending on $T>0$ only. Then, we obtain	\begin{align*}
	\mathbb{E}\left[\left\vert\frac{1}{n}\sum_{j=1}^n\check{U}_t^j-m_U^*(t)\right\vert^2 \right]&\leq \frac{1}{n}\sum_{j=1}^n\mathbb{E}\left[\left\vert\check{U}_t^j-m_U^*(t)\right\vert^2 \right]\leq \frac{2}{n}\sum_{j=1}^n\left\{\mathbb{E}\left[\sup_{t\in [0,T]}\vert\check{U}_t^j\vert^2 \right]+\|(m_U^*)^2\|_T \right\}\leq C.
\end{align*}
Thus, Lemma \ref{lemma2}-({\romannumeral1}) and the dominated convergence theorem (DCT) yield that $E_1^i\to 0$ as $n\to \infty$.        
Analogously, we can deduce that $E_3^i\to 0$ and $E_5^i\to 0$ as $n\to \infty$. For the term $E_2^i$, it holds that
    	\begin{align*}
    		E_2^i&=\beta\mathbb{E}\left[ \int_0^T\left\vert\left(\frac{1}{n}\sum_{j=1}^n\check{U}_t^j\right)^2-(m_U^*(t))^2\right\vert dt \right]\\
    		&\leq\beta\int_0^T\mathbb{E}\left[ \left\vert\left(\frac{1}{n}\sum_{j=1}^n\check{U}_t^j\right)^2-(\bar{U}_t^*)^2\right\vert \right]dt+\beta\int_0^T\mathbb{E}\left[\left\vert(\bar{U}_t^*)^2-(m_U^*(t))^2\right\vert \right]dt.
    	\end{align*}
Moreover,  using the inequality $(x-y)^2-x^2=y^2-2xy\leq y^2+2\vert x\vert \vert y \vert$, one has
\begin{align*}
&\mathbb{E}\left[\left\vert\left(\frac{1}{n}\sum_{j=1}^n\check{U}_t^j\right)^2-(\bar{U}_t^*)^2\right\vert\right]
\leq\mathbb{E}\left[\left\vert \bar{U}_t^*-\frac{1}{n}\sum_{j=1}^n\check{U}_t^j\right\vert^2\right]+ 2\left\{\mathbb{E}\left[\left\vert \bar{U}_t^*-\frac{1}{n}\sum_{j=1}^n\check{U}_t^j\right\vert^2\right]\mathbb{E}\left[\vert \bar{U}_t^*\vert^2\right]\right\}^{\frac{1}{2}}
\nonumber\\
&\quad\leq\mathbb{E}\left[\left\vert \bar{U}_t^*-\frac{1}{n}\sum_{j=1}^n\check{U}_t^j\right\vert^2\right]+ 2\left\{\mathbb{E}\left[\left\vert \bar{U}_t^*-\frac{1}{n}\sum_{j=1}^n\check{U}_t^j\right\vert^2\right]\frac{1}{n}\sum_{i=1}^n\mathbb{E}\left[\vert U_t^{*,i}\vert^2\right]\right\}^{\frac{1}{2}}\to 0, ~ \text{as}~ n\to \infty.
\end{align*}
Thus, employing the DCT again, we have $E_2^i\to 0$, as $n\to \infty$. The limit  $E_4^i\to 0$, as $n\to \infty$ can be verified similarly. Lastly, we obtain 
    	\begin{equation}\label{r1}
    		J_i(\theta^i, \boldsymbol{\theta}^{*, -i})-\bar{J}_i(\theta^i; m_U^*)\to 0, \   n\to \infty.
    	\end{equation}	
    	
Next, we consider the second term on the RHS of \eqref{A1}. To this purpose, it is necessary to introduce the following auxiliary stochastic control problem described as
\begin{equation}\label{problem2}
\left\{
\begin{aligned}
&\inf_{\theta^i\in\mathcal{A}}\tilde{J}_i(\theta^i; m_U^*):= \mathbb{E}\left[\int_0^T\left((\theta_t^i)^2+\rho\theta_t^i(\tilde{U}_t^i-m_U^*(t))+\beta\vert\tilde{U}_t^i-m_U^*(t)\vert^2\right)dt+\gamma\vert\tilde{U}_T^i-m_U^*(T)\vert^2\right],\\[0.4em]
&\tilde{U}_t^i=
u_i-a_i\int_0^t(\tilde{U}_s^i-m_U^*(s))ds+\int_0^tc_i\theta^i_sds-\int_0^t\int_0^1z\tilde{U}_{s-}^{i}N^i(ds, dz)+\int_0^tm_{\phi}^*(s)ds.
\end{aligned}\right.	
\end{equation}
Using a straight calculation as in Section \ref{sec:MFE}, we get that, the optimal feedback strategy of \eqref{problem2} is given by
    	\begin{equation*}
    		\tilde{\theta}^{*,i}_t=\theta^{*,i}(t, \tilde{U}_t^{*,i})
    		=\left(-c_iA_i(t)-\frac{\rho}{2}\right)(\tilde{U}_t^{*,i}-m_U^*(t))-\frac{c_i}{2}B_i(t, m_U^*, m_{\phi}^*),
    	\end{equation*}
where the state process $\tilde{U}^{*,i}=(\tilde{U}_t^{*,i})_{t\in[0,T]}$ satisfies the dynamics 
    	\begin{equation*}
    		\tilde{U}_t^{*,i}=u_i-a_i\int_0^t(\tilde{U}_s^{*,i}-m_U^*(s))ds+\int_0^tc_i\tilde{\theta}^{*,i}_sds-\int_0^t\int_0^1z\tilde{U}_{s-}^{*,i}N^i(ds, dz)+\int_0^tm_{\phi}^*(s)ds
    	\end{equation*}
with $A_i=(A_i(t))_{t\in[0,T]}$ and $B_i=(B_i(t))_{t\in[0,T]}$ being given in \eqref{npara1}. In other words, $\inf_{\theta^i\in\mathcal{A}}\tilde{J}_i(\theta^i; m_U^*)$ $=\tilde{J}_i(\tilde{\theta}^{*,i}; m_U^*)$. Then, we have 
\begin{align}\label{A2}
    	&\inf_{\theta^i\in\mathcal{A}}\bar{J}_i(\theta^i; m_U^*)-\bar{J}_i(\check{\theta}^{*,i}; m_U^*)
    		=\inf_{\theta^i\in\mathcal{A}}\bar{J}_i(\theta^i; m_U^*)-\inf_{\theta^i\in\mathcal{A}}\tilde{J}_i(\theta^i; m_U^*)+\inf_{\theta^i\in\mathcal{A}}\tilde{J}_i(\theta^i; m_U^*)-\bar{J}_i(\check{\theta}^{*,i}; m_U^*)\nonumber\\
    		&\qquad\qquad\qquad\geq \inf_{\theta^i\in\mathcal{A}}\left(\bar{J}_i(\theta^i; m_U^*)-\tilde{J}_i(\theta^i; m_U^*)\right)+\left(\tilde{J}_i(\tilde{\theta}^{*,i}; m_U^*)-\bar{J}_i(\check{\theta}^{*,i}; m_U^*)\right).
    	\end{align}
For the first term on the RHS of \eqref{A2}, we have
    	\begin{align*}
    		&\vert \bar{J}_i(\theta^{i},m_U^*)-\tilde{J}_i(\theta^{i},m_U^*) \vert \nonumber\\
            &\quad\leq \rho\mathbb{E}\left[\int_0^T\vert\theta_t^i(\check{U}_t^{i}-\tilde{U}_t^i)\vert dt \right]+\beta\mathbb{E}\left[\int_0^T\vert(\check{U}_t^{i})^2-(\tilde{U}_t^i)^2\vert dt\right]+2\beta\mathbb{E}\left[\int_0^T\vert m_U^*(t)(\check{U}_t^i-\tilde{U}_t^i)\vert dt\right]\nonumber\\
        &\qquad+\gamma\mathbb{E}\left[\vert(\check{U}_T^{i})^2-(\tilde{U}_T^i)^2\vert\right]+2\gamma\mathbb{E}\left[ \vert m_U^*(T)(\check{U}_T^i-\tilde{U}_T^i)\vert\right]=: I_1^i+I_2^i+I_3^i+I_4^i+I_5^i.
    	\end{align*}
It follows from \eqref{eq:boundednesstheta2} that $\mathbb{E}[\sup_{t\in[0,T]}\vert \tilde{U}_t^i\vert^2]\vee\mathbb{E}[\sup_{t\in[0,T]}\vert \check{U}_t^i\vert^2]\leq C$ for some finite $C>0$ depending on $T>0$ only. Then, we have from \eqref{checkU} and \eqref{problem2} that
    	\begin{align*}
    		&\mathbb{E}\left[\vert\check{U}_t^i-\tilde{U}_t^i\vert^2\right]\leq C\Bigg\{\int_0^t\mathbb{E}\left[\vert\check{U}_s^i-\tilde{U}_s^i\vert^2\right]ds+\int_0^t\mathbb{E}\left[\left\vert\frac{1}{n}\sum_{j=1}^n\check{U}_s^j-m_U^*(s)\right\vert^2\right]ds\\
            &\qquad+\mathbb{E}\left[\left\vert\frac{1}{n}\sum_{j\neq i}\int_0^t\int_0^1(kz\check{U}_s^j+\ell)N^j(ds,dz)-\int_0^t\left(\int_0^1z\nu(dz)km_U^*(s) +\ell\right)ds\right\vert^2\right]\Bigg\},
    	\end{align*}
where we have that
\begin{align*}
&\mathbb{E}\left[\left\vert\frac{1}{n}\sum_{j\neq i}\int_0^t\int_0^1(kz\check{U}_s^j+\ell)N^j(ds,dz)-\int_0^t\left(\int_0^1z\nu(dz)km_U^*(s)+\ell\right)ds \right\vert^2\right]\\
&\leq C\Bigg\{\mathbb{E}\left[\left\vert\frac{1}{n}\sum_{j\neq i}\int_0^t\int_0^1(kz\check{U}_s^j+\ell)N^j(ds,dz)-\frac{1}{n}\sum_{j\neq i}\int_0^t\int_0^1(kz\check{U}_s^j+\ell)\nu(dz)ds\right\vert^2\right]\\
&\quad+\mathbb{E}\left[\left\vert\frac{1}{n}\sum_{j\neq i}\int_0^t\int_0^1(kz\check{U}_s^j+\ell)\nu(dz)ds-\int_0^t\left(\int_0^1z\nu(dz)km_U^*(s)+\ell\right)ds \right\vert^2\right]\Bigg\}\\
&\leq C\left\{\frac{1}{n^2}+\frac{1}{n^2}\sum_{j\neq i}\int_0^t\mathbb{E}[\vert\check{U}_s^j\vert^2]ds
+\int_0^t\mathbb{E}\left[\left\vert\frac{1}{n}\sum_{j\neq i}(\check{U}_s^j-m_U^*(s))\right\vert^2\right]ds\right\}\\
&\leq C\left\{\frac{1}{n^2}+\frac{1}{n^2}\sum_{j=1}^n\int_0^t\mathbb{E}[\vert\check{U}_s^j\vert^2]ds
+\int_0^t\mathbb{E}\left[\left\vert\frac{1}{n}\sum_{j=1}^n\check{U}_s^j-m_U^*(s)\right\vert^2\right]ds\right\}\to 0,\ \text{as} \ n\to \infty.
\end{align*}
Using an analogous reasoning in Lemma  \ref{lemma2}-({\romannumeral3}) together with Gronwall’s lemma, we can claim that 
    	\begin{equation}\label{e3}
    		\lim_{n\to \infty}\mathbb{E}\left[\vert \check{U}_t^i-\tilde{U}_t^i\vert^2\right]= 0.
    	\end{equation}
For the term $I_1^i$, we have from \eqref{eq:boundednesstheta2} and DCT that
    	\begin{align*}
I_1^i&=\rho\mathbb{E}\left[\int_0^T\left\vert\theta_t^i(\check{U}_t^{i}-\tilde{U}_t^i)\right\vert dt \right]\leq \rho\left\{\mathbb{E}\left[\int_0^T\left\vert \check{U}_t^i-\tilde{U}_t^i\right\vert^2 dt\right]\left\|\theta^i\right\|_{2,T}^2\right\}^{\frac{1}{2}}\to 0, \ \text{as}\ n\to \infty.
\end{align*}
Similar to $I_1^i$, it is not difficult to deduce that $I_3^i\to 0$ and $I_5^i\to 0$ as $n\to\infty$. For the term $I_2^i$, it holds that
    	\begin{equation*}
I_2^i=\beta\int_0^T\mathbb{E}\left[\vert(\check{U}_t^{i})^2-(\tilde{U}_t^i)^2\vert \right]dt\leq \beta\int_0^T\left\{\mathbb{E}\left[ \vert \check{U}_t^{i}-\tilde{U}_t^i\vert^2 \right]+2\left\{\mathbb{E}\left[\vert \check{U}_t^i\vert^2 \right]\mathbb{E}\left[\vert \check{U}_t^{i}-\tilde{U}_t^i\vert^2 \right]\right\}^{\frac{1}{2}}\right\}dt.
\end{equation*}
Thus, $I_2^i\to 0$, as $n\to \infty$. Similarly, we have $I_4^i\to 0$, as $n\to \infty$. Consequently
\begin{align}\label{r2}
    \bar{J}_i(\theta^{i},m_U^*)-\tilde{J}_i(\theta^{i},m_U^*) \to 0, \ \text{as}  \ n\to \infty.
\end{align} 
    	For the second term on the RHS of \eqref{A2},  
    	\begin{align}\label{A3}
    		&\left\vert\tilde{J}_i(\tilde{\theta}^{*,i}; m_U^*)-\bar{J}_i(\check{\theta}^{*,i};m_U^*)\right\vert\nonumber\\
    		&=\Bigg\vert\mathbb{E}\Bigg[\int_0^T\Big((\tilde{\theta}_t^{*,i})^2-(\check{\theta}_t^{*,i})^2+\rho\tilde{\theta}_t^{*,i}(\tilde{U}_t^{*,i}-m_U^*(t))-\rho\check{\theta}_t^{*,i}(\check{U}_t^{*,i}-m_U^*(t))+\beta((\tilde{U}_t^{*,i})^2\nonumber\\
    		&\quad-(\check{U}_t^{*,i})^2)-2\beta m_U^*(t)(\tilde{U}_t^{*,i}-\check{U}_t^{*,i})\Big)dt+\gamma\left((\tilde{U}_T^{*,i})^2-(U_T^{*,i})^2-2 m_U^*(T)(\tilde{U}_T^{*,i}-U_T^{*,i})\right)\Bigg]\Bigg\vert\nonumber\\&\leq\mathbb{E}\left[\int_0^T\vert(\tilde{\theta}_t^{*,i})^2-(\check{\theta}_t^{*,i})^2\vert dt \right]+\rho\mathbb{E}\left[\int_0^T\vert\tilde{\theta}_t^{*,i}\tilde{U}_t^{*,i}-\check{\theta}_t^{*,i}\check{U}_t^{*,i}\vert dt \right]+\rho\mathbb{E}\left[\int_0^T\vert m_U^*(t)(\tilde{\theta}_t^{*,i}-\check{\theta}_t^{*,i}) \vert dt \right]\nonumber\\
    		&\quad +\beta\mathbb{E}\left[\int_0^T\vert (\tilde{U}_t^{*,i})^2-(\check{U}_t^{*,i})^2\vert dt \right]+2\beta\mathbb{E}\left[\int_0^T\vert m_U^*(t)(\tilde{U}_t^{*,i}-\check{U}_t^{*,i}) \vert dt \right]+\gamma\mathbb{E}\left[\vert(\tilde{U}_T^{*,i})^2-(U_T^{*,i})^2\vert\right]\nonumber\\
    		&\quad+2\gamma\mathbb{E}\left[\vert m_U^*(T)(\tilde{U}_T^{*,i}-U_T^{*,i})\vert\right]=:L_1^i+L_2^i+L_3^i+L_4^i+L_5^i+L_6^i+L_7^i.
    	\end{align}
 Recalled that $\check{\theta}^{*,i}_t
=(-c_iA_i(t)-\frac{\rho}{2})(\check{U}_t^{*,i}-m_U^*(t))-\frac{c_i}{2}B_i(t, m_U^*, m_{\phi}^*)$ and $\tilde{\theta}^{*,i}_t=(-c_iA_i(t)-\frac{\rho}{2})(\tilde{U}_t^{*,i}-m_U^*(t))-\frac{c_i}{2}B_i(t, m_U^*, m_{\phi}^*)$ for $t\in[0,T]$. Then, one has, for all $t\in[0,T]$,
\begin{equation}\label{l1}
    \mathbb{E}\left[\vert \tilde{\theta}_t^{*,i}-\check{\theta}_t^{*,i}\vert^2\right]\leq C \mathbb{E}\left[\vert \tilde{U}_t^{*,i}-\check{U}_t^{*,i}\vert^2\right].
\end{equation}
Using the argument which led to Lemma \ref{lemma2}-({\romannumeral2}), there exists a constant $C>0$ depending on $T>0$ only such that 
\begin{equation}\label{l2}
\mathbb{E}\left[\sup_{t\in [0,T]}\vert \tilde{U}_t^{*,i}\vert^2\right]\vee\mathbb{E}\left[\sup_{t\in [0,T]}\vert \tilde{\theta}_t^{*,i}\vert^2\right] \leq C.
\end{equation}
Following the proof procedure of \eqref{e3},  we obtain from Gronwall's lemma and the dominated convergence theorem that 
\begin{equation}\label{l3}
\mathbb{E}\left[\vert \tilde{U}_t^{*,i}-\check{U}_t^{*,i} \vert^2\right]\to 0, \ \text{as} \  n\to \infty.
\end{equation}
By applying the inequality $(x-y)^2-x^2=y^2-2xy\leq y^2+2\vert x\vert\vert y\vert$ and estimates \eqref{l1}-\eqref{l3}, we have
\begin{align*}
L_1^i &= \int_0^T \mathbb{E}\left[\left\vert(\tilde{\theta}_t^{*,i})^2-(\check{\theta}_t^{*,i})^2\right\vert\right] dt
\leq \int_0^T\left\{\mathbb{E}\left[\left\vert \tilde{\theta}_t^{*,i}-\check{\theta}_t^{*,i}\right\vert^2\right]+ 2\left(\mathbb{E}\left[\vert \tilde{\theta}_t^{*,i}-\check{\theta}_t^{*,i}\vert^2\right]\mathbb{E}\left[\vert  \tilde{\theta}_t^{*,i}\vert^2\right]\right)^{\frac{1}{2}} \right\}dt\nonumber\\
&\leq C\int_0^T\left\{\mathbb{E}\left[\vert \tilde{U}_t^{*,i}-\check{U}_t^{*,i}\vert^2\right]+ \left\{\mathbb{E}\left[\vert \tilde{U}_t^{*,i}-\check{U}_t^{*,i}\vert^2\right]\mathbb{E}\left[\vert  \tilde{\theta}_t^{*,i}\vert^2\right]\right\}^{\frac{1}{2}}\right\}dt\to 0, \ \text{as}\ n\to \infty.
\end{align*}
We then follow the similar argument in showing that $L_4^i \to 0$ and $L_6^i\to 0$ as $n\to\infty$. For the term $L_2^i$, we have
\begin{align*}
L_2^i&=\rho\int_0^T\mathbb{E}\left[\left\vert\tilde{\theta}_t^{*,i}\tilde{U}_t^{*,i}-\check{\theta}^{*,i}\check{U}_t^{*,i}\right\vert\right]dt=\rho\int_0^T\mathbb{E}\left[\left\vert\tilde{\theta}_t^{*,i}(\tilde{U}_t^{*,i}-\check{U}_t^{*,i})+(\tilde{\theta}_t^{*,i}-\check{\theta}_t^{*,i})\check{U}_t^{*,i}\right\vert\right]dt\nonumber\\
&\leq \rho\int_0^T\left\{\left\{\mathbb{E}\left[\vert\tilde{\theta}_t^{*,i}\vert^2\right]\mathbb{E}\left[\vert \tilde{U}_t^{*,i}-\check{U}_t^{*,i}\vert^2\right]\right\}^{\frac{1}{2}}+\left\{\mathbb{E}\left[\vert \tilde{\theta}_t^{*,i}-\check{\theta}_t^{*,i}\vert^2\right]\mathbb{E}\left[\vert \check{U}_t^{*,i}\vert^2\right]\right\}^{\frac{1}{2}}\right\}dt\nonumber\\
&\leq \rho\int_0^T\left\{\left\{\mathbb{E}\left[\vert\tilde{\theta}_t^{*,i}\vert^2\right]\mathbb{E}\left[\vert \tilde{U}_t^{*,i}-\check{U}_t^{*,i}\vert^2\right]\right\}^{\frac{1}{2}}+C\left\{\mathbb{E}\left[\vert \tilde{U}_t^{*,i}-\check{U}_t^{*,i} \vert^2\right]\mathbb{E}\left[\vert \check{U}_t^{*,i}\vert^2\right]\right\}^{\frac{1}{2}}\right\}dt\\
&\leq C\int_0^T\left\{\mathbb{E}\left[\vert \tilde{U}_t^{*,i}-\check{U}_t^{*,i} \vert^2\right]\right\}^{\frac{1}{2}}dt\to 0, \ \ \text{as} \ n\to 0.
\end{align*}
For the term $L_3^i$, we have 
\begin{align*}
L_3^i&=\rho\int_0^T
\mathbb{E}\left[\left\vert m_U^*(t)(\tilde{\theta}_t^{*,i}-\check{\theta}_t^{*,i})\right\vert\right]dt\leq \rho\int_0^T\left\{\mathbb{E}\left[\vert \tilde{\theta}_t^{*,i}-\check{\theta}_t^{*,i} \vert^2\right]\mathbb{E}\left[\vert m_U^*(t)\vert^2\right]\right\}^{\frac{1}{2}}dt\\
&\leq C\int_0^T\left\{\mathbb{E}\left[\left\vert \tilde{U}_t^{*,i}-\check{U}_t^{*,i} \right\vert^2\right]\mathbb{E}\left[\left\vert m_U^*(t)\right\vert^2\right]\right\}^{\frac{1}{2}}dt\to 0, \quad \text{as} \ n\to \infty.
\end{align*}
Similarly, we can obtain $L_5^i \to 0$, $L_7^i \to 0$, as $n\to \infty$. It follows that all terms on the RHS of \eqref{A3} converge to 0 as $n\to \infty$, which results in that
\begin{equation}\label{r3}
\tilde{J}_i(\tilde{\theta}^{*,i}; m_U^*)-\bar{J}_i(\check{\theta}^{*,i};m_U^*)\to 0, \quad \text{as} \ n\to \infty. 
\end{equation}
Lastly, we focus on the third term on the RHS of \eqref{A1}. In fact, we have
\begin{align*}
&\left\vert\bar{J}_i(\check{\theta}^{*,i}; m_U^*)-J_i(\boldsymbol{\theta}^*)\right\vert\nonumber\\
&=\Bigg\vert\mathbb{E}\Bigg[\int_0^T\Big((\check{\theta}_t^{*,i})^2-(\theta^{*,i}_t)^2+\rho\check{\theta}_t^{*,i}(\check{U}_t^{*,i}-m_U^*(t))-\rho\theta_t^{*,i}(U_t^{*,i}-m_U^*(t))+\beta((\check{U}_t^{*,i})^2\nonumber\\
&\quad-(U_t^{*,i})^2)-2\beta m_U^*(t)(\check{U}_t^{*,i}-U_t^{*,i})\Big)dt+\gamma\left((\check{U}_T^{*,i})^2-(U_T^{*,i})^2-2 m_U^*(T)(\check{U}_T^{*,i}-U_T^{*,i})\right)\Bigg]\Bigg\vert\nonumber\\
&\leq\mathbb{E}\left[\int_0^T\vert(\check{\theta}_t^{*,i})^2-(\theta_t^{*,i})^2\vert dt\right]+\rho\mathbb{E}\left[\int_0^T\vert\check{\theta}_t^{*,i}\check{U}_t^{*,i}-\theta_t^{*,i}U_t^{*,i}\vert dt\right]+\rho\mathbb{E}\left[\int_0^T\vert m_U^*(t)(\check{\theta}_t^{*,i}-\theta_t^{*,i}) \vert dt\right]\nonumber\\
&\quad +\beta\mathbb{E}\left[\int_0^T\vert (\check{U}_t^{*,i})^2-(U_t^{*,i})^2\vert dt\right]+2\beta\mathbb{E}\left[\int_0^T\vert m_U^*(t)(\check{U}_t^{*,i}-U_t^{*,i}) \vert dt\right]+\gamma\mathbb{E}\left[\vert(\check{U}_T^{*,i})^2-(U_T^{*,i})^2\vert\right]\nonumber\\
&\quad+2\gamma\mathbb{E}\left[\vert m_U^*(T)(\check{U}_T^{*,i}-U_T^{*,i})\vert\right].
\end{align*}
Using the results we obtained before, it can be demonstrated that 
\begin{equation}\label{r4}
\bar{J}_i(\check{\theta}^{*,i}; m_U^*)-J_i(\boldsymbol{\theta}^*)\to 0, \quad \text{as}\ n\to \infty.
    	\end{equation}
The proof of \eqref{r4} is identical to that of \eqref{r3}, and we omit it here. Thus, the limits \eqref{r1}, \eqref{r2}, \eqref{r3} and \eqref{r4} jointly imply \eqref{ANE}. Thus, we complete the proof of the theorem.
\end{proof}

	

\appendix
\section{Proofs of Auxiliary Results}\label{sec:appendix}   	

In this appendix, we collect the proof of the auxiliary Lemma \ref{lemma2}.
    
{\small\begin{proof}[Proof of Lemma \ref{lemma2}]
Throughout the proof, let $C$ be a positive constant depending on $T>0$ only, but which may be different from line to line.  (i) It follows from the dynamics \eqref{ndynamics} of the membrane potential process $U^i=(U_t^i)_{t\in[0,T]}$ that
	\begin{align*}
		\vert U_t^i \vert^2
		&\leq C\Bigg\{1+\int_0^t\vert\theta^i_s\vert^2ds+\int_0^t\vert U_s^i\vert^2ds+\int_0^t\frac{1}{n}\sum_{j=1}^n\vert U_s^j\vert^2ds+\left\vert\int_0^t\int_0^1zU_{s-}^i\tilde{N}^i(ds,dz)\right\vert^2\nonumber\\
		&\quad+\left\vert\frac{1}{n}\sum_{j\neq i}\int_0^t\int_0^1 \phi( zU_{s-}^j)\tilde{N}^j(ds,dz)\right\vert^2\Bigg\},
	\end{align*}
	where $\tilde{N}^i(ds,dz):=N^i(ds,dz)-\nu(dz)ds$ is the compensate Poisson random measure. Taking expectation on both sides of the above display and applying BDG inequality, we have
	\begin{align}\label{moment2}
	&\mathbb{E}\left[\sup_{t\in [0, T]}\vert U_t^i\vert^2\right]
    \leq C\left\{ 1+ \left\|\theta^i\right\|_{2,T}^2 +\int_0^T \mathbb{E}\left[\sup_{s\in [0, t]}\vert U_s^i\vert^2\right] dt+\int_0^T\frac{1}{n}\sum_{j=1}^n\mathbb{E}\left[\sup_{s\in [0, t]}\vert U_s^j \vert^2\right]dt\right\}.
	\end{align}
	Therefore, it yields that
	\begin{equation*}
		\frac{1}{n}\sum_{i=1}^n \mathbb{E}\left[\sup_{t\in [0, T]}\vert U_t^i \vert^2\right]\leq C\left\{1+\frac{1}{n}\sum_{i=1}^n\left\|\theta^i\right\|_{2,T}^2+\int_0^T\frac{1}{n}\sum_{j=1}^n \mathbb{E}\left[	\sup_{s\in [0, t]}\vert U_s^j \vert^2\right]dt \right\}.
	\end{equation*}
	Thus,  by the Gronwall's lemma with the $L^2$-boundedness of strategy $\theta^i$ in our assumption, we have that $\frac{1}{n}\sum_{i=1}^n\mathbb{E}[\sup_{t\in [0, T]}\vert U_t^i \vert^2]\leq Ce^{CT}<\infty.$ Plugging this estimate into \eqref{moment2}, we obtain that $\mathbb{E}\left[\sup_{t\in [0, T]} \vert U_t^i\vert^2\right] $ $\leq C\left\{1+\int_0^T	\mathbb{E}\left[\sup_{s\in [0, t]}\vert U_s^i\vert^2\right]dt\right\}$. Then, the desired estimate of item (i) follows from the Gronwall's lemma again. 
	
	For item (ii), it follows from the strategy $\theta^{*,i}$ which is defined by \eqref{ntheta} that $\mathbb{E}\left[\vert \theta^{*,i}_t \vert^2\right]\leq C\{1+\mathbb{E}\left[\vert U_t^{*,i}\vert^2\right]\}$ and 
	\begin{align*}
		\mathbb{E}\left[\sup_{t\in [0, T]}\vert U_t^{*,i} \vert^2\right]&\leq C\left\{1+\int_0^T\mathbb{E}\left[\sup_{s\in [0, t]}\vert U_s^{*,i}\vert^2\right]dt+\int_0^T\frac{1}{n}\sum_{j=1}^n\mathbb{E}\left[\sup_{s\in [0, t]}\vert U_s^{*,j}\vert^2\right]dt+\left\|\theta^i\right\|_{2,T}^2\right\}.
	\end{align*}
	Then, the desired result follows from the similar argument as item (i). 
	
	Lastly, the proof of item (iii) is more complex. Let us recall that  $m_U^*\in\mathcal{C}_T^1$ satisfies that
	\begin{align*}
		m_{U}^*(t)
		&=\mathbb{E}_{\mu}[\boldsymbol{u}]-\frac{1}{2}\int_0^t\mathbb{E}_{\mu}[\boldsymbol{c}^2B_{\boldsymbol{p}}(s,m_U^{*, (\boldsymbol{p})},m_{\phi}^{*, (\boldsymbol{p})})]ds-\int_0^1z\nu(dz)\int_0^tm_U^*(s)ds+\int_0^tm_{\phi}^*(s)ds.
	\end{align*}
	We introduce an auxiliary process $(\hat{U}_t^{*,i})_{t\in[0,T]}$ as follows, for $i=1,\ldots, n$, 
	\begin{align*}
		\hat{U}_t^{*,i}&=\mathbb{E}_{\mu}[\boldsymbol{u}]-\int_0^t\mathbb{E}_{\mu}\left[\boldsymbol{a}+\boldsymbol{c}^2A_{\boldsymbol{p}}(s)+\frac{\rho \boldsymbol{c}}{2}\right](\hat{U}_s^{*,i}-m_U^*(s))ds-\frac{1}{2}\int_0^t\mathbb{E}[\boldsymbol{c}^2B_{\boldsymbol{p}}(s,m_U^{*,(\boldsymbol{p})},m_{\phi}^{*,(\boldsymbol{p})})]ds\nonumber\\&\quad-\int_0^t\int_0^1z\hat{U}_{s-}^{*,i}N^i(ds, dz)+\int_0^tm_{\phi}^*(s)ds,
	\end{align*}
	and an alternative auxiliary process $(\tilde{U}_t^{*,i})_{t\in[0,T]}$ given by, for $i=1,\ldots, n$, 
	\begin{align*}
		\tilde{U}_t^{*,i}&=u_i-\int_0^t\left(a_i+c_i^2A_i(s)+\frac{\rho c_i}{2}\right)(\tilde{U}_s^{*,i}-m_U^*(s))ds-\frac{c_i^2}{2}\int_0^tB_i(s,m_U^*,m_{\phi}^*)ds\nonumber\\&\quad-\int_0^t\int_0^1z\tilde{U}_{s-}^{*, i}N^i(ds, dz)+\int_0^tm_{\phi}^*(s)ds.
	\end{align*}
	Then, it follows from the Gronwall's lemma together with the boundedness of $m_U^*$ and $m_{\phi}^*$ on $[0, T]$ that $\mathbb{E}\left[\sup_{t\in[0,T]}\vert \hat{U}_t^{*,i} \vert^2\right]\leq C$. And we can deduce that
	\begin{align*}
		\mathbb{E}\left[\sup_{t\in[0,T]}\vert \tilde{U}_t^{*,i} \vert^2\right]\leq C\left\{1+\int_0^T\mathbb{E}\left[\sup_{s\in[0, t]}\vert \tilde{U}_s^{*,i} \vert^2\right]dt +\int_0^T\frac{1}{n}\sum_{j=1}^n\mathbb{E}\left[\sup_{s\in[0,t]}\vert \tilde{U}_s^{*,j} \vert^2\right]dt\right\}.
	\end{align*}
	This implies from the Gronwall's lemma that $\mathbb{E}[\sup_{t\in[0,T]}\vert \tilde{U}_t^{*,i} \vert^2]\leq C$. On the other hand, we can make a decomposition as follows:
	\begin{align}\label{main}
		&\mathbb{E}\left[\left\vert\bar{U}_t^*-m_U^*(t)\right\vert^2\right]=\mathbb{E}\left[\left\vert \frac{1}{n}\sum_{i=1}^n(U_t^{*,i}-\tilde{U}_t^{*,i})+\frac{1}{n}\sum_{i=1}^n(\tilde{U}_t^{*,i}-\hat{U}_t^{*,i})+\frac{1}{n}\sum_{i=1}^n(\hat{U}_t^{*,i}-m_U^*(t))\right\vert^2\right]\nonumber\\
		&\leq 3\mathbb{E}\left[\left\vert \frac{1}{n}\sum_{i=1}^n(U_t^{*,i}-\tilde{U}_t^{*,i})\right\vert^2\right]+3\mathbb{E}\left[\left\vert\frac{1}{n}\sum_{i=1}^n(\tilde{U}_t^{*,i}-\hat{U}_t^{*,i})\right\vert^2\right]+3\mathbb{E}\left[\left\vert\frac{1}{n}\sum_{i=1}^n(\hat{U}_t^{*,i}-m_U^*(t))\right\vert^2\right].
	\end{align}
	Consider the first term on the RHS of \eqref{main}, we can deduce that
	\begin{align*}
		\mathbb{E}\left[\left\vert \frac{1}{n}\sum_{i=1}^n(U_t^{*,i}-\tilde{U}_t^{*,i})\right\vert^2\right] &\leq C\left\{ \int_0^T 	\mathbb{E}\left[\left\vert\frac{1}{n}\sum_{i=1}^n( U_s^{*,i}-\tilde{U}_s^{*,i}) \right\vert^2\right]ds+\int_0^T\mathbb{E}\left[\vert \bar{U}^*_s-m_U^*(s)\vert^2\right]ds\right.\nonumber\\&\quad+\mathbb{E}\left[\left\vert\frac{1}{n}\sum_{i=1}^n\int_0^t\int_0^1z(U_{s-}^{*,i}-\tilde{U}_{s-}^{*,i})N^i(ds, dz)\right\vert^2\right]\nonumber\\&\left.\quad+\mathbb{E}\left[\left\vert\frac{n-1}{n^2}\sum_{j=1}^n\int_0^t\int_0^1\phi(zU_{s-}^{*,j})N^j(ds,dz)-\int_0^tm_{\phi}^*(s)ds\right\vert^2\right]\right\}.
	\end{align*}
	As a result, using the isometry relation, it yields that
	\begin{align}\label{estimate}
		&\mathbb{E}\left[\left\vert\frac{1}{n}\sum_{i=1}^n( U_t^{*,i}-\tilde{U}_t^{*,i})\right\vert^2\right]\leq C\left\{\frac{1}{n}+\int_0^t\mathbb{E}\left[\left\vert \frac{1}{n}\sum_{i=1}^n(U_s^{*,i}-\tilde{U}_s^{*,i}) \right\vert^2\right]ds+\int_0^t\mathbb{E}\left[\left\vert\bar{U}_s^*-m_U^*(s) \right\vert^2\right]ds\right\}\nonumber\\
		&\leq C\left\{\frac{1}{n}+\int_0^t\mathbb{E}\left[\left\vert\frac{1}{n}\sum_{j=1 }^n( U_{s}^{*,j}-\tilde{U}_s^{*,j})\right\vert^2\right]ds+\int_0^t\mathbb{E}\left[\left\vert\frac{1}{n}\sum_{j=1}^n(\tilde{U}_s^{*,j}-\hat{U}_s^{*,j})\right\vert^2\right]ds\right.\nonumber\\
		&\qquad\quad\left.+\int_0^t\mathbb{E}\left[\left\vert\frac{1}{n}\sum_{j=1}^n(\hat{U}_s^{*,j}-m_U^*(s))\right\vert^2\right]ds\right\}.
	\end{align}
	For the second term on the RHS of \eqref{main}, we can deduce that
	\begin{align*}
		&\tilde{U}_t^{*,i}-\hat{U}_t^{*,i}=u_i-\mathbb{E}_{\mu}[\boldsymbol{u}]-\int_0^t\left(a_i+c_i^2A_i(s)+\frac{\rho c_i}{2}\right)\tilde{U}_s^{*,i}ds+\int_0^t\mathbb{E}_{\mu}\left[\boldsymbol{a}+\boldsymbol{c}^2A_{\boldsymbol{p}}(s)+\frac{\rho \boldsymbol{c}}{2}\right]\hat{U}_s^{*,i}ds\nonumber\\
		&\quad-\frac{1}{2}\int_0^tc_i^2B_i(s, m_U^*,m_{\phi}^*)ds+\frac{1}{2}\int_0^t\mathbb{E}_{\mu}[\boldsymbol{c}^2B_{\boldsymbol{p}}(s,m_U^{*,(\boldsymbol{p})},m_{\phi}^{*,(\boldsymbol{p})})]ds-\int_0^t\int_0^1z(\tilde{U}_{s-}^{*,i}-\hat{U}_{s-}^{*,i})N^i(ds, dz).
	\end{align*}
	Then, we have
	\begin{align}\label{estimate2}
		&\mathbb{E}\left[\left\vert\frac{1}{n}\sum_{i=1}^n(\tilde{U}_t^{*,i}-\hat{U}_t^{*,i})\right\vert^2\right]\nonumber\\&\leq C\Bigg\{\left\vert\frac{1}{n}\sum_{i=1}^n(u_i-\mathbb{E}_{\mu}[\boldsymbol{u}])\right\vert^2+\left\vert\frac{1}{n}\sum_{i=1}^n\int_0^tc_i^2B_i(s, m_U^*, m_{\phi}^*)ds-\int_0^t\mathbb{E}_{\mu}[\boldsymbol{c}^2B_{\boldsymbol{p}}(s,m_U^{*,(\boldsymbol{p})},m_{\phi}^{*,(\boldsymbol{p})})]ds\right\vert^2\nonumber\\&\quad+\mathbb{E}\left[\left\vert\frac{1}{n}\sum_{i=1}^n\int_0^t\left(a_i+c_i^2A_i(s)+\frac{\rho c_i}{2}\right)\tilde{U}_s^{*,i}ds-\frac{1}{n}\sum_{i=1}^n\int_0^t\mathbb{E}_{\mu}\left[\boldsymbol{a}+\boldsymbol{c}^2A_{\boldsymbol{p}}(s)+\frac{\rho \boldsymbol{c}}{2}\right]\hat{U}_s^{*,i}ds\right\vert^2\right]\nonumber\\&\quad+\mathbb{E}\left[\left\vert\frac{1}{n}\sum_{i=1}^n\int_0^t\int_0^1z(\tilde{U}_{s-}^{*,i}-\hat{U}_{s-}^{*,i})N^i(ds, dz)\right\vert^2\right]\Bigg\}.
	\end{align}
	For the first term on the RHS of \eqref{estimate2}, when the assumption \ref{Ap} is in force, we have 
	\begin{equation}\label{t1}
		\left\vert\frac{1}{n}\sum_{i=1}^n u_i-\mathbb{E}_{\mu}[\boldsymbol{u}] \right\vert^2=\left\vert \int_{\mathcal{O}}u\mu_n(dp)-\int_{\mathcal{O}}u\mu(dp)\right\vert^2\to 0, \quad \text{as}\ n\to \infty.
	\end{equation}
	For the second term on the RHS of \eqref{estimate2}, we have
    \begin{align}
    	&\left\vert\frac{1}{n}\sum_{i=1}^n\int_0^tc_i^2B_i(s, m_U^*, m_{\phi}^*)ds-\int_0^t\mathbb{E}_{\mu}[\boldsymbol{c}^2B_{\boldsymbol{p}}(s,m_U^{*,(\boldsymbol{p})},m_{\phi}^{*,(\boldsymbol{p})})]ds\right\vert^2\nonumber\\&\leq
    	2\int_0^t\left\vert\int_{\mathcal{O}}\left(\frac{1}{n}\sum_{i=1}^nc_i^2B_i(s, m_U^{*,(p)}, m_{\phi}^{*,(p)})-c^2B_i(s,m_U^{*,(p)},m_{\phi}^{*,(p)})\right)\mu(dp)\right\vert^2ds\nonumber\\&\quad +2\int_0^t\left\vert \int_{\mathcal{O}}c^2\left(\frac{1}{n}\sum_{i=1}^nB_i(s, m_U^{*,(p)}, m_{\phi}^{*,(p)})-B_p(s,m_U^{*,(p)},m_{\phi}^{*,(p)})\right)\mu(dp)\right\vert^2ds\nonumber\\&\leq 
    	C\left\{\int_0^t\left\vert\frac{1}{n}\sum_{i=1}^nc_i^2-\int_{\mathcal{O}}c^2\mu(dp)\right\vert^2ds+\int_0^t\left\vert \frac{1}{n}\sum_{i=1}^nh_i(s)-\int_{\mathcal{O}}h_p(s)\mu(dp)\right\vert^2ds\right.\nonumber\\&\quad\left.+\int_0^t\int_s^T\left| \frac{1}{n}\sum_{i=1}^n\frac{A_i(v)}{h_i(v)}-\int_{\mathcal{O}}\frac{A_p(v)}{h_p(v)}\mu(dp)\right|^2dvds\right\}.
    \end{align}
    By using Assumption \ref{Ap}, and note that $\mu_n$ is the empicical measure for $n\geq1$, it follows from DCT that, as $n\to\infty$, 
    \begin{equation}
    	\int_0^t\left\vert\frac{1}{n}\sum_{i=1}^nc_i^2-\int_{\mathcal{O}}c^2\mu(dp)\right\vert^2ds=\int_0^t\left\vert\int_{\mathcal{O}}c^2\mu_n(dp)-\int_{\mathcal{O}}c^2\mu(dp)\right\vert^2ds\to 0.
    \end{equation}
    Then, we shall verify that for fixed $t\in [0, T]$, $A_p(t)\in C_b(\mathcal{O})$ and $h_p(t)\in C_b(\mathcal{O})$. Indeed, for a given sequence $(p_{_k})_{k\geq1}\subset\mathcal{O}$ converging to some $p\in\mathcal{O}$ in the sense of $\|p_{_k}-p\|_1:=|u_{k}-u|+|a_{k}-a|+|c_{k}-c|\to 0$ as $k \rightarrow \infty$, it holds that  for any $t\in [0, T]$,
	\begin{align*}
		&\left\vert A_{p_{_k}}(t)-A_p(t)\right\vert\leq 
		\left\vert \frac{\delta_{p_{_k}}^+}{c_{_k}^2}-\frac{\delta_p^+}{c^2}\right\vert+\left\vert\left( \frac{\delta_{p_{_k}}^+}{c_{_k}^2}-\frac{\delta_p^+}{c^2}\right)\frac{\delta_p^+-\delta_p^-}{(c^2\gamma-\delta_p^-)e^{(\delta_p^+-\delta_p^-)(T-t)}-(c^2\gamma-\delta_p^+)}\right\vert\\
		&\qquad +\left\vert \frac{\delta_{p_{_k}}^+}{c_{_k}^2}\left(\frac{\delta_{p_{_k}}^+-\delta_{p_{_k}}^-}{(c_{_k}^2\gamma-\delta_{p_{_k}}^-)e^{(\delta_{p_{_k}}^+-\delta_{p_{_k}}^-)(T-t)}-(c_{_k}^2\gamma-\delta_{p_{_k}}^+)}-\frac{\delta_p^+-\delta_p^-}{(c^2\gamma-\delta_p^-)e^{(\delta_p^+-\delta_p^-)(T-t)}-(c^2\gamma-\delta_p^+)}\right)\right\vert\\
		&\quad\leq C\left(\vert \delta_{p_{_k}}^+-\delta_p^+\vert+\vert \delta_{p_{_k}}^--\delta_p^-\vert+\vert c_k-c\vert\right)\leq C\Vert p_{_k}-p \Vert_1.
	\end{align*}
	It is also not difficult to obtain $\vert h_{p_{_k}}(t)-h_p(t)\vert\to 0$ as $\|p_k-p\|_1\to 0$. And furthermore, we have $A$ and $h$ are continuous on the compact set $[0,T]\times\mathcal{O}$.  It follows from DCT that, as $n\to\infty$, 
	\begin{align}
	&	\int_0^t\left\vert \frac{1}{n}\sum_{i=1}^nh_i(s)-\int_{\mathcal{O}}h_p(s)\mu(dp)\right\vert^2ds=\int_0^t\left\vert \int_{\mathcal{O}}h_p(s)\mu_n(dp)-\int_{\mathcal{O}}h_p(s)\mu(dp)\right\vert^2ds\to 0,\\
	&	\int_0^t\int_s^T\left| \frac{1}{n}\sum_{i=1}^n\frac{A_i(v)}{h_i(v)}-\int_{\mathcal{O}}\frac{A_p(v)}{h_p(v)}\mu(dp)\right|^2dvds=\int_0^t\int_s^T\left|\int_{\mathcal{O}}\frac{A_p(v)}{h_p(v)}\mu_n(dp) -\int_{\mathcal{O}}\frac{A_p(v)}{h_p(v)}\mu(dp)\right|^2dvds\to 0.\nonumber
	\end{align}
	For the third term on the RHS of \eqref{estimate2}, it can be deduced that
	\begin{align}\label{t3}
		&\mathbb{E}\left[\left\vert\frac{1}{n}\sum_{i=1}^n\int_0^t\left(a_i+c_i^2A_i(s)+\frac{\rho c_i}{2}\right)\tilde{U}_s^{*,i}ds-\int_0^t\mathbb{E}_{\mu}\left[\boldsymbol{a}+{\boldsymbol{c}}^2A_{\boldsymbol{p}}(s)+\frac{\rho \boldsymbol{c}}{2}\right]\hat{U}_s^{*,i}ds\right\vert^2\right]\nonumber\\ 
        &= \mathbb{E}\Bigg[\Bigg\vert\frac{1}{n}\sum_{i=1}^n\int_0^t\left(a_i+c_i^2A_i(s)+\frac{\rho c_i}{2}\right)(\tilde{U}_s^{*,i}-\hat{U}_s^{*,i})ds+\frac{1}{n}\sum_{i=1}^n\int_0^t\Big(\left(a_i+c_i^2A_i(s)+\frac{\rho c_i}{2}\right)\nonumber\\
        &\qquad-\mathbb{E}_{\mu}\left[\boldsymbol{a}+{\boldsymbol{c}}^2A_{\boldsymbol{p}}(s)+\frac{\rho \boldsymbol{c}}{2}\right]\Big)\hat{U}_s^{*,i}ds\Bigg\vert^2\Bigg]\\ 
		&\leq C\left\{\int_0^t\left\vert\frac{1}{n}\sum_{i=1}^n\left(a_i+c_i^2A_i(s)+\frac{\rho c_i}{2}\right)-\mathbb{E}_{\mu}\left[\boldsymbol{a}+{\boldsymbol{c}}^2A_{\boldsymbol{p}}(s)+\frac{\rho \boldsymbol{c}}{2}\right]\right\vert^2ds+\int_0^t\mathbb{E}\left[\left\vert \frac{1}{n}\sum_{i=1}^n(\tilde{U}_s^{*,i}-\hat{U}_s^{*,i})\right\vert^2\right]ds\right\}.\nonumber 
	\end{align}
	Set $G_1(t,p)=a+c^2A_p(t)+\frac{\rho c}{2}$ for $p\in \mathcal{O}$. Similarly, we can obtain that $G_1(\cdot,\cdot)$ is continuous on the compact set  $[0,T]\times\mathcal{O}$. Then, we have from DCT that, as $n\to\infty$,
	\begin{align}\label{t4}
		&\int_0^t\left\vert\frac{1}{n}\sum_{i=1}^n\left(a_i+c_i^2A_i(s)+\frac{\rho c_i}{2}\right)-\mathbb{E}_{\mu}\left[\boldsymbol{a}+{\boldsymbol{c}}^2A_{\boldsymbol{p}}(s)+\frac{\rho \boldsymbol{c}}{2}\right]\right\vert^2ds\\
		&\quad=\int_0^t\left\vert\int_{\mathcal{O}}G_1(s,p)\mu_n(dp)-\int_{\mathcal{O}}G_1(s,p)\mu(dp)\right\vert^2ds
		\to 0.\nonumber
	\end{align}
	For the last term on the RHS of \eqref{estimate2}, by the BDG inequality, it gives that
	\begin{equation}\label{t5}
		\mathbb{E}\left[\left\vert\frac{1}{n}\sum_{i=1}^n\int_0^t\int_0^1z(\tilde{U}_{s-}^{*,i}-\hat{U}_{s-}^{*,i})N^i(ds, dz)\right\vert^2\right]\leq C \int_0^t \mathbb{E}\left[\left\vert\frac{1}{n}\sum_{i=1}^n(\tilde{U}_{s-}^{*,i}-\hat{U}_{s-}^{*,i})\right\vert^2\right]ds.
	\end{equation}
	Plugging \eqref{t1}-\eqref{t5} into \eqref{estimate2}, it reduces to that
	\begin{align*}
		\mathbb{E}\left[\left\vert\frac{1}{n}\sum_{i=1}^n(\tilde{U}_t^{*,i}-\hat{U}_t^{*,i})\right\vert^2\right]&\leq C\int_0^t\mathbb{E}\left[\left\vert\frac{1}{n}\sum_{i=1}^n(\tilde{U}_s^{*,i}-\hat{U}_s^{*,i})\right\vert^2\right]ds.
	\end{align*}
	Utilizing the Gronwall's lemma, we have
	\begin{equation}\label{term2}
		\mathbb{E}\left[\left\vert\frac{1}{n}\sum_{i=1}^n(\tilde{U}_t^{*,i}-\hat{U}_t^{*,i})\right\vert^2\right]\to 0, \quad \text{as}\ \ n\to \infty.
	\end{equation}
	
	Now  we focus on the third term on the RHS of \eqref{main}, for fixed $t \in [0,T]$,  $(\hat{U}_t^{*,i})_{i=1}^n $ are i.i.d. copies with mean $m_U^*$. Thus, for $t\in[0,T]$,
	\begin{equation}\label{term3}
		\mathbb{E}\left[\left\vert\frac{1}{n}\sum_{i=1}^n(\hat{U}_t^{*,i}-m_U^*(t))\right\vert^2 \right]=\frac{1}{n^2}\sum_{i=1}^n\mathbb{E}\left[\vert\hat{U}_t^{*,i}-m_U^*(t)\vert^2 \right]=\frac{1}{n}\mathbb{E}\left[\vert\hat{U}_t^{*,1}-m_U^*(t)\vert^2 \right]\to 0, \ \text{as}\ n\to \infty.
	\end{equation}
	Plugging \eqref{term2} and \eqref{term3} into \eqref{estimate}, it becomes that
	\begin{align*}
		\mathbb{E}\left[\left\vert \frac{1}{n}\sum_{i=1}^n(U_t^{*,i}-\tilde{U}_t^{*,i})\right\vert^2\right]\leq C\int_0^t\mathbb{E}\left[\left\vert \frac{1}{n}\sum_{i=1}^n(U_t^{*,i}-\tilde{U}_t^{*,i})\right\vert^2\right]ds \rightarrow 0,   \ \text{as}\ n \rightarrow \infty.
	\end{align*} 
	Employing the Gronwall's lemma again, we have for any  $t \in[0,T]$,
	\begin{equation}\label{term1}
		\mathbb{E}\left[\left\vert \frac{1}{n}\sum_{i=1}^n(U_t^{*,i}-\tilde{U}_t^{*,i})\right\vert^2\right]\to 0, \ \text{as}\  n\to \infty.
	\end{equation}
	Finally, \eqref{term2}, \eqref{term3} and \eqref{term1} jointly yield that $\mathbb{E}\left[\vert\bar{U}_t^*-m_U^*(t)\vert^2\right]\to 0$ as $n\to\infty$.
	Thus, we complete the proof of the lemma.
\end{proof}
	}

\end{document}